\documentclass[times]{nlaauth}
\usepackage{amsmath,amssymb,graphics,epsfig,color,cite}
\usepackage[ruled,vlined]{algorithm2e}

\newtheorem{theorem}{Theorem}[section]
\newtheorem{proposition}[theorem]{Proposition}

\theoremstyle{remark}

\theoremstyle{definition}

\numberwithin{equation}{section}
\numberwithin{theorem}{section}

\newcommand{\mc}[1]{{\mathcal #1}}
\newcommand{\bb}[1]{{\mathbb #1}}

\newcommand{\rme}{\mathrm{e}}
\newcommand{\rmi}{\mathrm{i}}

\newcommand{\Lameps}{\Lambda_\varepsilon}
\begin{document}
\title{Approximated structured pseudospectra}
\runningheads{S. Noschese and L. Reichel}{Approximated structured pseudospectra}

\author{Silvia Noschese\affil{1}\corrauth and Lothar Reichel\affil{2}}

\address{\affilnum{1}Dipartimento di Matematica, SAPIENZA Universit\`a di Roma,
P.le Aldo Moro 5, 00185 Roma, Italy.\break
\affilnum{2}Department of Mathematical Sciences, Kent State University,
Kent, OH 44242, USA.}

\corraddr{Dipartimento di Matematica, SAPIENZA Universit\`a di Roma,
P.le Aldo Moro 5, 00185 Roma, Italy. E-mail: noschese@mat.uniroma1.it}



\begin{abstract}
Pseudospectra and structured pseudospectra are important tools for the analysis of 
matrices. Their computation, however, can be very demanding for all but small matrices.
A new approach to compute approximations of pseudospectra and structured pseudospectra,
based on determining the spectra of many suitably chosen rank-one or projected rank-one 
perturbations of the given matrix is proposed. The choice of rank-one or projected 
rank-one perturbations is inspired by Wilkinson's analysis of eigenvalue sensitivity. 
Numerical examples illustrate that the proposed approach gives much better insight into 
the pseudospectra and structured pseudospectra than random or structured random rank-one 
perturbations with lower computational burden. The latter approach is presently commonly 
used for the determination of structured pseudospectra.
\end{abstract}

\maketitle

\section{Introduction}\label{sec:1}
Many applications in science and engineering require knowledge of the location of some
or all eigenvalues of a matrix and the sensitivity of the eigenvalues to perturbations of
the matrix. The sensitivity can be studied with the aid of the eigenvalue condition 
number, based on particular rank-one perturbations of the matrix, as described by 
Wilkinson \cite[Chapter 2]{W65}, or by computing pseudospectra. Let $\Lambda(A)$ denote 
the spectrum of the matrix $A\in\mathbb{C}^{n\times n}$. The $\varepsilon$-pseudospectrum 
of the matrix $A$ is defined as
\begin{equation}\label{pseu}
\Lameps(A):=\left\{\lambda\in\mathbb{C}\colon\lambda\in\Lambda(A+E),~~
E\in\mathbb{C}^{n\times n},~~\|E\|\le\varepsilon\right\}
\end{equation}
for some $\varepsilon>0$. An insightful discussion of the $\varepsilon$-pseudospectrum and
many applications are presented by Trefethen and Embree \cite{TE}. The matrix norm 
$\|\cdot\|$ in \eqref{pseu} often is chosen to be the spectral norm. However, it will be 
convenient to instead use the Frobenius norm $\|\cdot\|_F$ in the present paper. Thus, for 
$E=[e_{ij}]_{i,j=1}^n\in\mathbb{C}^{n\times n}$, we have 
$\|E\|_F=\sqrt{\sum_{i,j=1}^n|e_{ij}|^2}$. The set \eqref{pseu} depends on the choice of 
matrix norm, however, this dependence often is not important in applications when one is 
interested in determining which eigenvalues of the matrix $A$ are most sensitive to 
perturbations. 

Algorithms for eigenvalue computations that respect the matrix structure may yield higher
accuracy and require less computing time than structure-ignoring methods. They also may
preserve eigenvalue symmetries in finite-precision arithmetic. A structure-respecting 
eigenvalue algorithm is said to be strongly backward stable, if the computed eigenvalues 
are exact eigenvalues of a slightly perturbed matrix with the same structure as the 
original matrix; see, e.g., Bunch \cite{B87}. To assess the numerical properties of a 
structure-respecting eigenvalue algorithm suitable measures of the sensitivity of the 
eigenvalues should be used in order not to overestimate the worst-case effect of 
perturbations. These measures include structured condition numbers, see Higham and Higham 
\cite{HH} and Karow et al. \cite{KKT}, as well as the structured 
$\varepsilon$-pseudospectrum. The latter can be applied to measure the sensitivity of the 
eigenvalues of a structured matrix to similarly structured perturbations. It is defined as
follows. Let ${\mc S}$ denote the subset of matrices in $\mathbb{C}^{n\times n}$ with a 
particular structure, such as bandedness, Toeplitz, Hankel, or Hamiltonian. Then, for some
$\varepsilon>0$, the structured $\varepsilon$-pseudospectrum of a matrix $A\in{\mc S}$ is 
given by
\begin{equation}\label{strctpseu}
\Lameps^{\mc S}(A):=\left\{\lambda\in\mathbb{C}\colon\lambda\in\Lambda(A+E),~~
E\in{\mc S},~~\|E\|\le\varepsilon\right\};
\end{equation}
see, e.g., \cite{BGK01,G06,R06} for discussions and illustrations. 

The computation of the (standard) $\varepsilon$-pseudospectrum \eqref{pseu} for a large or 
moderately sized matrix $A\in\mathbb{C}^{n\times n}$ for a fixed $\varepsilon>0$ can be 
very time-consuming. For instance, when the norm $\|\cdot\|$ in \eqref{pseu} is the 
spectral norm, approximations of the $\varepsilon$-pseudospectrum often are determined by 
computing the smallest singular value of many matrices of the form $A-z I_n$, where 
$I_n\in\mathbb{C}^{n \times n}$ denotes the identity matrix and $z\in\mathbb{C}$. If the
smallest singular value is smaller than or equal to $\varepsilon$, then $z$ belongs to the
set \eqref{pseu}. These computations are very demanding unless $A$ is small. To reduce the 
computational burden somewhat, it is suggested in \cite{WT,TE} that one first computes the
Schur factorization $A = U R U^H$ and then determines the smallest singular value of the 
matrix $R-zI_n$ for many $z$-values in $\mathbb{C}$. Here $U\in\mathbb{C}^{n\times n}$ is a
unitary matrix, $R\in\mathbb{C}^{n\times n}$ an upper triangular matrix, and the 
superscript $^H$ denotes transposition and complex conjugation. Nevertheless, the 
computational task is substantial also when applying the Schur factorization of a 
moderately sized or large matrix $A$. Moreover, the Schur factorization of $A$ cannot be
applied for the computation of the structured $\varepsilon$-pseudospectrum 
\eqref{strctpseu}. In fact, there are few methods available for computing the structured 
$\varepsilon$-pseudospectrum besides plotting the spectra of matrices 
$A\in\mathbb{C}^{n\times n}$ with structured random perturbations. Similarly, there are 
few methods available for determining the stability radius  under structured 
perturbations. 
In fact, the computation of structured 
$\varepsilon$-pseudospectra has become an established tool in gaining 
insight into behavior of matrix-based models in dynamical system theory 
under structured perturbations.
We recall that the  (structured) stability radius is the smallest value of 
$\varepsilon$ for which a (structured) $\varepsilon$-pseudospectrum contour reaches the 
imaginary axis and defines the norm of the smallest (structured) perturbation that destroys 
the  (structured) stability, that is, having all the eigenvalues confined to $\mathbb{C}^{-}$.
Just as the spectral abscissa of a matrix provides a measure of its stability, that is, the asymptotic 
decay of associated dynamical systems, so does the (structured) $\varepsilon$-pseudospectral abscissa,
i.e. the maximal real part  of  points of the (structured) $\varepsilon$-pseudospectrum
provides a measure of robust (structured) stability, where by robust we mean with respect to 
(structured) perturbations of the matrix.
The computations of these quantities remain significant computational 
challenges to date. 

The high computational burden of computing standard (unstructured) pseudospectra has 
spurred the development of algorithms that can be executed efficiently on a parallel
computer; see, e.g., Bekas and Gallopoulos \cite{BG} and Mezher and Philippe \cite{MP}.
We propose a different approach to speed up the computations, that also can be applied to
the determination of structured pseudospectra. The given matrix 
$A\in\mathbb{C}^{n\times n}$ is modified by a sequence of rank-one matrices that are known
to yield large perturbations of the eigenvalues according to Wilkinson's analysis 
\cite[Chapter 2]{W65}. Already fairly few rank-one matrices provide insight into the
$\varepsilon$-pseudospectrum and when different components of the pseudospectrum coalesce.
Computed examples illustrate that the number of our chosen rank-one matrices required to
gain knowledge of the $\varepsilon$-pseudospectrum is much smaller than when random 
rank-one perturbations are used. Our method can be used to inexpensively compute 
approximated standard (unstructured) pseudospectra when available software tools, such
as Eigtool \cite{Wr02} and Seigtool \cite{KKK10}, are too expensive to use. 

To determine approximations of structured pseudospectra, we project the rank-one matrices 
suggested by Wilkinson's analysis onto the set of matrices with desired structure. The 
computations with our approach can be implemented efficiently on parallel computers, but 
we will not pursue this aspect in the present paper.

We remark that  
Karow \cite{K10} has analyzed structured pseudospectra when the underlying norm is 
unitarily invariant. The distance to standard and structured defectivity is not 
considered in \cite{K10}. Also the computational approach of the present paper is new.   
Rump \cite{R06} has characterized the structured pseudospectrum for some 
structures, such as Toeplitz, and has by means of computer-assisted proofs shown that 
the structured $\varepsilon$-pseudospectrum \eqref{strctpseu} is quite similar to the 
standard $\varepsilon$-pseudospectrum \eqref{pseu} when there is no zero-structure such 
as bandedness.  Rump's analysis does not apply to banded Toeplitz matrices. A linearly 
convergent algorithm  for the computation of the structured $\varepsilon$-pseudospectral abscissa and 
radius of a Toeplitz matrix, and for determining sections of $\Lameps^{\mc S}(A)$ near 
extremal points is described in \cite{BGN}. Computations with this algorithm generally 
are time-consuming. 

This paper is organized as follows. Section \ref{sec1b} discusses rank-one perturbations 
and presents some background material. In Section \ref{sec:2}, we provide formulas for 
eigenvalue structured condition numbers and for the maximal structured perturbations for 
Toeplitz and Hamiltonian structures. These results are applied in Section \ref{sec:3} in 
the design of our approach for the inexpensive computation of approximated unstructured 
and structured pseudospectra. Numerical examples are presented in Section \ref{sec:4}, and
conclusions are provided in Section \ref{sec:5}.

\section{Rank-one perturbations and structured matrices}\label{sec1b}
The points in a structured $\varepsilon$-pseudospectrum \eqref{strctpseu} are exact eigenvalues 
of a nearby matrix in ${\mc S}$. This suggests that we may use standard results from the 
literature on eigenvalue sensitivity to infinitely small structured perturbations. The 
structured condition number of an eigenvalue $\lambda$ of $A$ is a first-order measure of 
the worst-case effect on $\lambda$ of perturbations with the same structure as $A$. The 
structured condition numbers used in this paper can easily be computed when endowing the
subspace of matrices with the Frobenius norm.

First consider an unstructured matrix $A\in\mathbb{C}^{n\times n}$ and assume that it has
the simple eigenvalue $\lambda$ with unit right and left eigenvectors $x$ and $y$, 
respectively. Let $E\in\mathbb{C}^{n\times n}$ have norm $\|E\|_F=1$ and assume that 
$\varepsilon>0$ is small enough so that the eigenvalue $\lambda_E(t)$ of $A+tE$ exists and
is unique for all $0\leq t<\varepsilon$. Then
\[
\lambda_E(t)=\lambda+\frac{y^H E x}{y^H x}t + {\mathcal O}(t^2);
\]
see Wilkinson \cite[Chapter 2]{W65} for details. We have 
\[
\left|\frac{y^H E x}{y^H x}\right|\leq\frac{1}{|y^H x|}
\]
with equality for 
\begin{equation}\label{rank1}
E:=\eta\, yx^H
\end{equation}
for any unimodular $\eta\in\mathbb{C}$. We refer to matrices of the form \eqref{rank1} as
Wilkinson perturbations and to the set
\begin{equation}\label{Wdisk}
{\mc D}(\lambda,t):=\left\{z\in\mathbb{C}: |z-\lambda|\leq \frac{1}{|y^H x|}t\right\}
\end{equation}
as the Wilkinson disk associated with $\lambda$ of radius $t\geq 0$. The condition number 
of the eigenvalue $\lambda$ is defined as
\begin{equation}\label{condnumb}
\kappa(\lambda):=\frac{1}{|y^H x|}.
\end{equation}

We turn to structured matrices in a set 
${\mc S}\,{\scriptscriptstyle{\substack{\subset \\ \ne}}}\,\mathbb{C}^{n \times n}$ and 
consider structured perturbations in ${\mc S}$. Let $M|_{\mc S}$ denote the matrix in
$\mc S$ closest to $M\in\mathbb{C}^{n\times n}$ with respect to the Frobenius norm. This
projection is used in the numerator of the eigenvalue condition number for structured 
perturbations; see \cite{NP06,NP07} and Section \ref{sec:2}. In particular, the condition
number for structured perturbations is smaller than the condition number for unstructured
perturbations. We also will use the normalized projection,
\[
M|_{\widehat{\mc S}} := \frac{M|_\mc S}{\|M|_\mc S\|_F}\;.
\]

In this paper, we are mainly concerned with banded or general Toeplitz structure, 
${\mc S}:={\mc T}$, and Hamiltonian structure, ${\mc S}:={\mc H}$. Toeplitz matrices arise 
in many applications including the solution of ordinary differential equations. It is 
interesting to investigate the sensitivity of the eigenvalues of a Toeplitz matrix with 
respect to finite structure-preserving perturbations and, in particular, the sensitivity 
of the rightmost eigenvalue. In \cite{BGN}, the authors computed the rightmost points of 
the structured pseudospectrum of a Toeplitz matrix and investigated the structured 
pseudospectrum  of tridiagonal Toeplitz matrices. The eigenvalues and eigenvectors of a
tridiagonal Toeplitz matrix are known in closed form and many quantities required for the 
analysis are easily computable \cite{NPR}. Moreover, the $\varepsilon$-pseudospectrum of a 
tridiagonal Toeplitz matrix is well approximated by ellipses as $\varepsilon$ approaches 
zero and the order $n$ of the matrix goes to infinity \cite{RT}.

The structure ${\mc T}$ is determined by the location of the nonzero diagonals of the
Toeplitz matrix. Since the points in a structured pseudospectrum are eigenvalues of a 
nearby structured matrix with the same zero diagonals as $A\in{\mc T}$, it is 
straightforward to verify that the matrix $M|_\mc T$ for an arbitrary matrix 
$M\in\mathbb{C}^{n\times n}$ is obtained by replacing all elements in a nonzero diagonal 
of $M$ by their arithmetic mean \cite{NP07}.

This construction of the closest matrix in ${\mc T}$ to a given matrix 
$M\in\mathbb{C}^{n\times n}$ can be generalized to Hankel matrices by considering 
anti-diagonals in place of diagonals. Several other structures, such as persymmetry and
skew-persymmetry, can be handled similarly; see \cite{NP07} for illustrations.

We turn to Hamiltonian structure. Let ${\mc H}$ denote the linear subspace of Hamiltonian 
matrices of order $2n$, i.e., 
\[
\begin{split}
\mc  H & := \left\{Q \in \bb C^{2n\times 2n} \colon QJ=(QJ)^H\right\} \\ & =
\left\{Q=\begin{pmatrix} K & M \\ L & -K^H \end{pmatrix} \colon K,L,M 
\in \bb C^{n\times n}\;,\;L=L^H\;,\;M=M^H\right\}\,,
\end{split}
\]
where 
\begin{equation}\label{Jmat}
J:=\begin{pmatrix} 0 & I_n \\ -I_n & 0 \end{pmatrix}
\end{equation}
is the fundamental symplectic matrix.
Hamiltonian eigenvalue problems arise from applications in systems and control theory; 
see, e.g., \cite{BKM05} and references therein. The sensitivity of the eigenvalues of a 
Hamiltonian matrix with respect to a Hamiltonian perturbation was studied in \cite{BN} and
computable formulas for the structured condition numbers were derived. An expression for 
the closest Hamiltonian matrix $M|_\mc H$ to a given matrix $M\in \mathbb{C}^{n \times n}$ 
is shown in the following section.

\section{Structured maximal perturbations of Toeplitz and Hamiltonian eigenvalue problems}
\label{sec:2}
The following proposition summarizes results from \cite{NP07} for Toeplitz matrices and
perturbations, and will be used in the sequel. 

\begin{proposition}\label{lem:p1}
Let $\lambda$ be a simple eigenvalue of a Toeplitz matrix 
$A\in{\mc T}\subset\mathbb{C}^{n\times n}$ with right and left eigenvectors $x$ and $y$, 
respectively, normalized so that $\|x\|_F = \|y\|_F=1$. Given any matrix 
$E\in{\mc T}$ with $\|E\|_F = 1$, let $\lambda_E(t)$ be an eigenvalue of $A+tE$ converging
to $\lambda$ as $t\to 0$. Then
\[
|\dot\lambda_E(0)| \le \max\left\{\left|\frac{y^HEx}{y^Hx}\right| , \; \|E\|_F = 1,\, 
E \in \mc T\,\right\} = \frac{\|yx^H|_{\mc T}\|_F}{|y^Hx|}
\]
and
\begin{eqnarray}\label{f2}
\nonumber
\dot\lambda_E(0) & = & \frac{\|yx^H|_{\mc T}\|_F}{|y^Hx|} \qquad \mathrm{if} \qquad 
E=\eta yx^H|_{\widehat{\mc T}},
\end{eqnarray}
for any unimodular $\eta \in \mathbb(C)$. Here $\dot\lambda_E(t)$ denotes the derivative of $\lambda_E(t)$ with respect to the
parameter $t$.
\end{proposition}
It follows from Proposition \ref{lem:p1} that the Toeplitz structured condition number is
given by
\[
\kappa^{\mc  T}(\lambda) := \frac{\|yx^H|_{\mc  T}\|_F}{|y^Hx|}\,.
\]

We turn to Hamiltonian matrices and perturbations.

\begin{proposition}[\cite{BN}]\label{prop:1}
The closest Hamiltonian matrix to a given matrix 
\[
A = \begin{pmatrix} A_1 & A_3 \\ A_2 & A_4 \end{pmatrix}\in \bb C^{2n\times 2n}
\]
with respect to the Frobenius norm is 
\[
A|_{\mc  H} = \frac 12 \begin{pmatrix} A_1-A_4^H & A_3+A_3^H \\ 
A_2+A_2^H & A_4-A_1^H \end{pmatrix} = \frac 12 (A+JA^H J),
\]
where $J$ is the fundamental symplectic matrix \eqref{Jmat} and 
$A_k\in\mathbb{C}^{n\times n}$ for $1\leq k\leq 4$.
\end{proposition}

\medskip

\begin{proposition}[\cite{BN}]\label{teo:1}
Let $\lambda$ be a simple eigenvalue of a Hamiltonian matrix 
$Q\in{\mc H}\subset\mathbb{C}^{2n\times 2n}$ with right and left eigenvectors $x$ and $y$,
respectively, normalized so that
\begin{equation}\label{p:1}
\|x\|_F=\|y\|_F = 1\,, \qquad {\rm Im}(y^HJx) = 0,
\end{equation}
where $J$ is defined by \eqref{Jmat}. Given any matrix $E\in{\mc H}$ with $\|E\|_F = 1$,
let $\lambda_E(t)$ be an eigenvalue of $Q+tE$ converging to $\lambda$ as $t\to 0$. Then
\begin{equation}\label{p:2-}
|\dot\lambda_E(0)| \le \max\left\{\left|\frac{y^HEx}{y^Hx}\right| \colon  \|E\|_F = 1,\, 
E\in\mc  H\,\right\} = \frac{\|yx^H|_{\mc  H}\|_F}{|y^Hx|}\,.
\end{equation}
Moreover, 
\begin{equation}
\label{p:2}
\dot\lambda_E(0) = \frac{\|yx^H|_{\mc  H}\|_F}{y^Hx} \qquad \mathrm{if} \qquad 
E=\pm yx^H|_{\widehat{\mc H}}\,.
\end{equation}
\end{proposition}

We obtain from Proposition \ref{teo:1} that the Hamiltonian structured condition number 
is given by
\[
\kappa^{\mc  H}(\lambda) := \frac{\|yx^H|_{\mc  H}\|_F}{|y^Hx|}\,.
\]
We remark that the worst-case effect perturbation turns out to be a rank-2 complex matrix. 
Moreover, since
the (unstructured) condition number of a simple eigenvalue is
$\kappa(\lambda) = \|yx^H\|_F/|y^Hx|$, we have 
$\kappa^{\mc  H}(\lambda) = \kappa(\lambda)$ if $yx^H$ is Hamiltonian. This can occur only
if $\lambda$ is a purely imaginary eigenvalue; see \cite{BN}.

\section{Approximated structured $\varepsilon$-pseudospectra}\label{sec:3}
This section describes how the structured Wilkinson perturbations of Section \ref{sec:2} 
can be applied to determine useful approximations of $\Lameps^{\mc S}(A)$ when $A$ is a
matrix in
${\mc S}$ 
with all eigenvalues distinct.
When ${\mc S}=\mathbb{C}^{n \times n}$, i.e., when 
$A$ has no particular structure, the perturbations that affect the eigenvalue $\lambda$ of
$A$ the most, relative to the norm of the perturbation, are multiples of the rank-one 
matrices \eqref{rank1}. The Wilkinson disks \eqref{Wdisk} for the different eigenvalues 
are disjoint if the radius $t$ of the disks is smaller than the distance $\varepsilon_{*}$
from defectivity of the matrix $A$,
\[
\varepsilon_{*} = \inf \{ \| A - B \|_F \colon B \in \mathbb{C}^{n\times n} \ \mbox{is defective} \}\;.
\]
Analogously, in case 
${\mc S}\,{\scriptscriptstyle{\substack{\subset \\ \ne}}}\,\mathbb{C}^{n \times n}$ the
threshold is the structured distance from defectivity $\varepsilon_*^{{\mc S}}$ of $A$,
\[
\varepsilon_*^{{\mc S}} =  \inf \{ \| A - B \|_F \colon B \in \mc S \ \mbox{is defective} \}\;.
\]
Clearly, $\varepsilon_*^{{\mc S}}\geq\varepsilon_{*}$; see, e.g., \cite{Dem83,AB,ABBO11,BGMN} 
for details. In the structured case, the rank-one Wilkinson perturbations \eqref{rank1}
are projected as described in Section \ref{sec:2}.

Assume that machine epsilon, $\varepsilon_M$, satisfies $0<\varepsilon_M\ll\varepsilon_{*}$.
First let ${\mc S}=\mathbb{C}^{n\times n}$. Then the component of 
$\Lambda_{\varepsilon_M}(A)$ that contains $\lambda $ is approximately a disk of radius
$\kappa(\lambda)\varepsilon_M=\varepsilon_M/|y^Hx|$ centered at $\lambda$, i.e. the Wilkinson disk  
${\mc D}(\lambda,\varepsilon_M)$. An estimate of
$\varepsilon_{*}$ is given by
\begin{equation}\label{rout_1}
\varepsilon:=\min_{\substack{1\le i\le n\\ 1\le j\le n \\ j\neq i}}
\frac{|\lambda_i -\lambda_j|}{\kappa(\lambda_i)+\kappa(\lambda_j)}\,.
\end{equation}
Indeed, ${\mc D}(\lambda_i,\varepsilon)$ is tangential to ${\mc D}(\lambda_j,\varepsilon)$ when 
$|\lambda_i -\lambda_j|=(\kappa(\lambda_i)+\kappa(\lambda_j))\,\varepsilon$. Let 
the index pair $\{\hat\imath,\hat\jmath\}$ minimize the ratio \eqref{rout_1} over all 
distinct eigenvalue pairs. Then the Wilkinson disks ${\mc D}(\lambda_{\hat\imath},t)$ and 
${\mc D}(\lambda_{\hat\jmath},t)$ are the disks that will coalesce first when 
increasing $t$. The union of the Wilkinson disks 
${\mc D}(\lambda_{\hat\imath},\varepsilon)$ and 
${\mc D}(\lambda_{\hat\jmath},\varepsilon)$ determine a rough approximation of the 
$\varepsilon$-pseudospectrum around the eigenvalues $\lambda_{\hat\imath}$ and 
$\lambda_{\hat\jmath}$ for $\varepsilon$ sufficiently small. We will refer to the 
eigenvalues $\lambda_{\hat\imath}$ and $\lambda_{\hat\jmath}$ as the 
{\it most $\Lameps$-sensitive pair of eigenvalues}. We note that usually the most 
$\Lameps$-sensitive pair of eigenvalues are not the two worst conditioned ones. 

We turn to the situation when 
${\mc S}\,{\scriptscriptstyle{\substack{\subset \\ \ne}}}\,\mathbb{C}^{n \times n}$. Then
the role of $\kappa(\lambda)$ is played by the first-order measure in the Frobenius norm 
of the worst-case effect on $\lambda$ of structured perturbations, 
i.e., the structured condition
number $\kappa^{\mc S}(\lambda)$. We refer to the set
\begin{equation*}\label{WSdisk}
{\mc D}^{{\mc S}}(\lambda,t):=\left\{z\in\mathbb{C}: |z-\lambda|\leq \frac{\|yx^H|_{\mc  S}\|_F}{|y^Hx|}t\right\}
\end{equation*}
as the ${\mc S}$-structured Wilkinson disk associated with $\lambda$ of radius $t\geq 0$.
For $\varepsilon_M\ll\varepsilon_{*}^{{\mc S}}$, the component of 
$\Lambda_{\varepsilon_M}^{\mc S}(A)$ that contains $\lambda $ is approximately a disk of
radius $\kappa^{\mc S}(\lambda)\varepsilon_M=
\kappa(\lambda)\|yx^H|_{\mc S}\|_F\,\varepsilon_M$ centered at $\lambda$,  i.e.,  the ${\mc S}$-structured Wilkinson disk  
${\mc D}^{{\mc S}}(\lambda,\varepsilon_M)$, and an estimate 
of $\varepsilon_{*}^{{\mc S}}$ is given by
\begin{equation}\label{rout_2}
\varepsilon^{\mc S}:= \min_{\substack{1\le i\le n\\ 1\le j\le n \\ j\neq i}}
\frac{|\lambda_i -\lambda_j|}{\kappa^{\mc S}(\lambda_i)+\kappa^{\mc S}(\lambda_j)}\geq
\varepsilon \,,
\end{equation}
and ${\mc D}^{{\mc S}}(\lambda_i,\varepsilon^{{\mc S}})$ is tangential to ${\mc D}^{{\mc S}}(\lambda_j,\varepsilon^{{\mc S}})$ when 
$|\lambda_i -\lambda_j|=(\kappa^{\mc S}(\lambda_i)+\kappa^{\mc S}(\lambda_j))\,\varepsilon^{{\mc S}}$. Let $\{\hat\imath,\hat\jmath\}$ be a minimizing index pair over all distinct pairs 
of eigenvalues. Then,  the ${\mc S}$-structured Wilkinson disks ${\mc D}^{{\mc S}}(\lambda_{\hat\imath},t)$ and 
${\mc D}^{{\mc S}}(\lambda_{\hat\jmath},t)$ are the first ones to coalesce as $t$ increases.
In the sequel, we will refer to the eigenvalues $\lambda_{\hat\imath}$ and 
$\lambda_{\hat\jmath}$ as the {\it most $\Lameps^{\mc S}$-sensitive pair of eigenvalues}; they are
not necessarily the worst conditioned eigenvalues with respect to structure-preserving 
perturbations $E\in {\mc S}$.

We found in numerous computations that the rank-one perturbation $E$ with all elements equal
to $\varepsilon/n$ generally induces a significant perturbation in the rightmost 
eigenvalue and gives a meaningful lower bound for the $\varepsilon$-pseudospectral abscissa. In the 
structured case, such a rank-one perturbation has to be projected as discussed above in 
order to give a useful approximation of the structured $\varepsilon$-pseudospectral abscissa. This is 
illustrated in Example 3 of the following section.

\section{Numerical illustrations}\label{sec:4}
This section presents computations that illustrate the approaches for determining 
approximated pseudospectra discussed in the previous section. All computations were 
carried out in MATLAB with about $16$ significant decimal digits. Throughout this
section $\rmi=\sqrt{-1}$.

\subsection{Toeplitz structure}
We consider two banded Toeplitz matrices, one of which is tridiagonal. There are well-known
explicit formulas for the eigenvalues and eigenvectors for the latter kind of matrix; see, 
e.g., \cite{NPR}. 

\begin{table}[htb!]
\centering
\begin{tabular}{cccc}\hline 
$i$ &$\lambda_i$ & $\kappa(\lambda_i)$ & $\kappa^{\mc T}(\lambda_i)$  \\ 
\hline  
$1$ &$-0.4988$ & $1.153\cdot 10^{2}$ & $2.625\cdot 10^{0}$ \\
$2$ &${\phantom 0}0.0564$ & $3.269\cdot 10^{2}$ & $1.559\cdot 10^{0}$ \\
$3$ &${\phantom 0}0.8147$ & $4.243\cdot 10^{2}$ & ${\phantom 0}4.472\cdot 10^{-1}$ \\
$4$ &${\phantom 0}1.5731$ & $3.269\cdot 10^{2}$ & $1.559\cdot 10^{0}$ \\
$5$ &${\phantom 0}2.1283$ & $1.153\cdot 10^{2}$ & $2.625\cdot 10^{0}$ \\
\hline
\end{tabular}
\caption{Example 1: Eigenvalue condition numbers.}
\label{Tab1}
\end{table}

{\bf Example 1}.
Consider a real tridiagonal Toeplitz matrix of order $n=5$ with random diagonal and 
superdiagonal entries in the interval $[0,1]$, and random subdiagonal entries in the 
interval $[0,5]$. This gives a matrix with fairly ill-conditioned eigenvalues. It is 
shown in \cite{NPR} that the sensitivity of the eigenvalues grows exponentially with the
ratio of the absolute values of the sub- and super-diagonal matrix entries. 

The eigenvalues for a typical tridiagonal Toeplitz matrix of the kind described and their 
standard and structured condition numbers are shown in Table \ref{Tab1}. While the 
eigenvalues in the middle of the spectrum are the worst conditioned with respect to 
unstructured perturbations, the extremal eigenvalues are most sensitive to structured 
perturbations. This is also discussed in \cite{NPR}.

\begin{figure}[ht]
\begin{center}
\includegraphics[width=6.5cm]{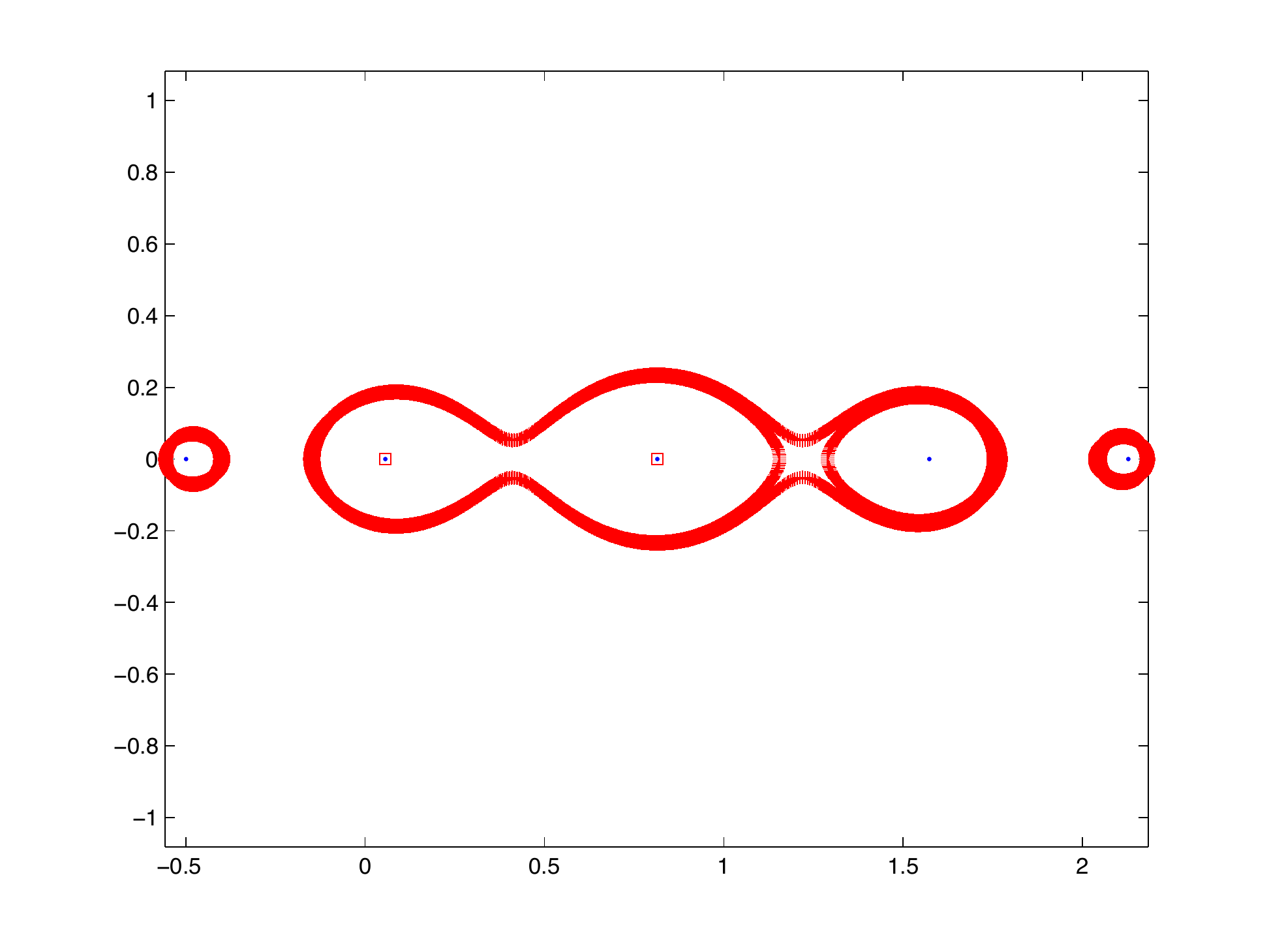}
\hskip 10mm
\includegraphics[width=6.5cm]{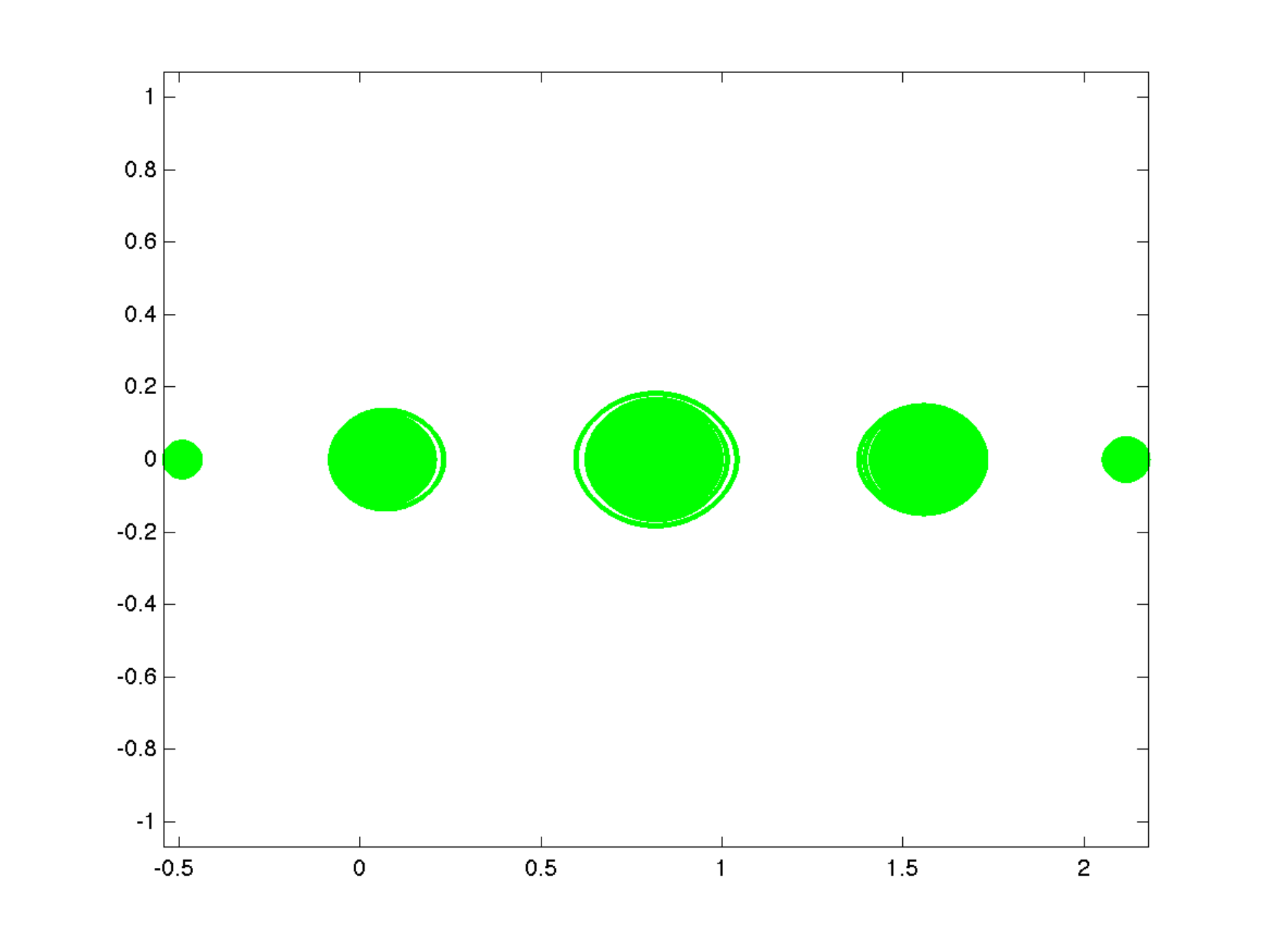}
\end{center}
\caption{Example 1. Left plot: $\Lambda_{\varepsilon_1}(A)$ is approximated by the 
eigenvalues of matrices of the form $A+\varepsilon_1 W_2$ and $A+\varepsilon_1 W_3$, where 
the $W_j$ are Wilkinson perturbations associated with the eigenvalues $\lambda_j$, $j=2,3$ 
(marked by red squares), for $\eta:=\rme^{\rmi\theta_k}$, $\theta_k:=2\pi(k-1)/10^3$, 
$k=1:10^3$, and $\varepsilon_1=10^{-3.2}$. Right plot: $\Lambda_{\varepsilon_1}(A)$ 
approximated by the eigenvalues of matrices of the form 
$A+\varepsilon_1\rme^{\rmi\theta_k} E_i$, $i,k=1:10^3$, where the $E_i$ are unit-norm 
rank-one random perturbations.\hfill\break}\label{fig1}
\end{figure}

\begin{figure}[ht]
\begin{center}
\includegraphics[width=13cm]{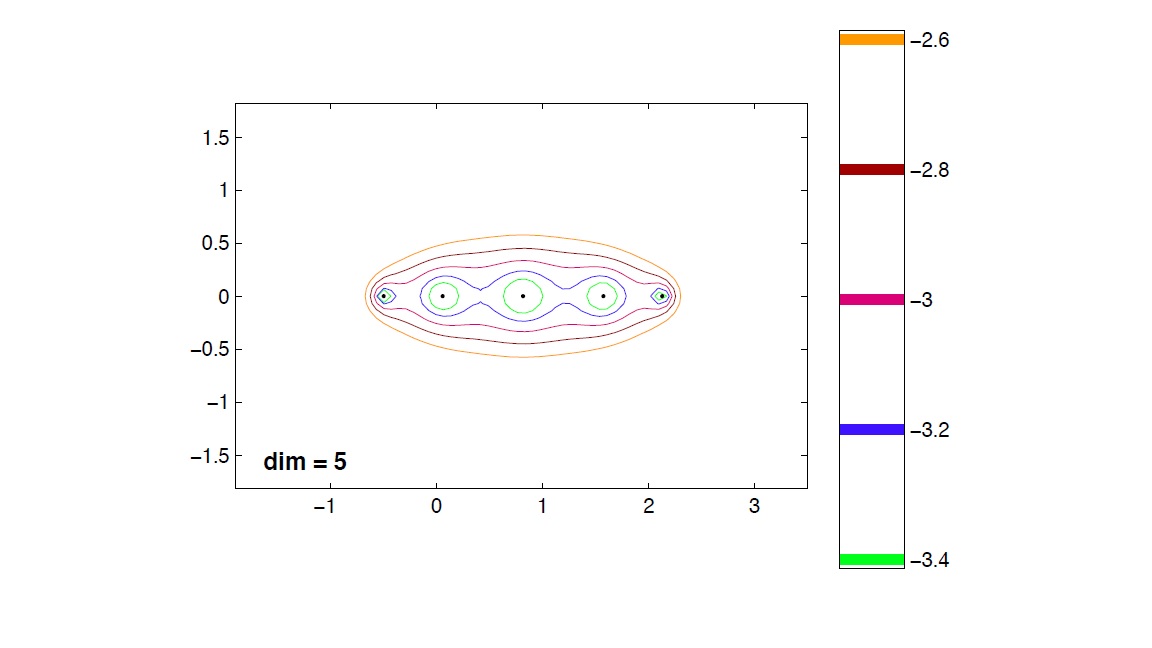}
\end{center}
\caption{Example 1: $\varepsilon$-pseudospectra by Eigtool, where $\varepsilon=10^k$, 
$k=-3.4:0.2:-2.6$. dim shows the order of the matrix.\hfill\break}\label{fig2}
\end{figure}

The estimate (\ref{rout_1}) of the (unstructured) distance from defectivity 
$\varepsilon_*$ is $\varepsilon_1=10^{-3.2}$. It is achieved for the indices $2$ and $3$ 
of the most $\Lameps$-sensitive pair of eigenvalues. Figure \ref{fig1} (left) displays the 
spectrum of matrices of the form $A+\varepsilon_1 W_2$ and $A+\varepsilon_1 W_3$, where
the $W_j$ are Wilkinson perturbations \eqref{rank1} associated with the eigenvalues 
$\lambda_j$, $j=2,3$, for $\eta:=\rme^{\rmi\theta_k}$ and $\theta_k:=2\pi(k-1)/10^3$, 
$k=1:10^3$. Here and throughout this section $\eta$ is the leading coefficient of the
Wilkinson perturbation \eqref{rank1}. Details of the computations are described by
Algorithm 1.

\begin{algorithm}
\DontPrintSemicolon
\KwData{matrix $A$, eigensystem $\{\lambda_i, x_i, y_i,\; \forall i=1:n\}$}
\KwResult{approximated $\Lambda_{\varepsilon}(A)$}

\nl compute $\varepsilon$, $\{\hat\imath,\hat\jmath\}$  by \eqref{rout_1}\;
\nl compute  $W_{\hat\imath}:=\varepsilon \,y_{\hat\imath}x_{\hat\imath}^H$, $W_{\hat\jmath}:=\varepsilon \,y_{\hat\jmath}x_{\hat\jmath}^H$  \;
\nl display the spectrum of  $A+\eta W_{\hat\imath}$ for $\eta:=\rme^{\rmi\theta_k}$, $\theta_k=2\pi(k-1)/10^3$, $k=1:10^3$ \;
\nl display the spectrum of  $A+\eta W_{\hat\jmath}$ for $\eta:=\rme^{\rmi\theta_k}$, $\theta_k=2\pi(k-1)/10^3$, $k=1:10^3$ \;
\caption{Algorithm for computing an approximated pseudospectrum} \label{algo1}
\end{algorithm}

Figure \ref{fig1} (right) displays the approximated
$\varepsilon_1$-pseudospectrum given by the spectra of matrices of the form 
$A+\varepsilon_1\rme^{\rmi\theta_k}E_i$, $i,k=1:10^3$, where the $\theta_k$ are
defined as above and the $E_i$ are random rank-one perturbations with
$\|E_i\|_F=1$. Thus, the figure shows spectra of $10^6$ matrices.  Figure \ref{fig2}
displays pseudospectra determined by Eigtool \cite{Wr02}. Comparing the 
$\varepsilon_1$-pseudospectrum of Figure \ref{fig2} with Figure \ref{fig1} illustrates the 
effectiveness of the simple approach proposed in this paper. In particular, the 
approximated $\varepsilon_1$-pseudospectrum of Figure \ref{fig1} (left) provides a much
better approximation of the $\varepsilon_1$-pseudospectrum than the approximated 
$\varepsilon_1$-pseudospectrum of Figure \ref{fig1} (right) and requires the computation
of many fewer spectra ($10^3$ versus $10^6$). This makes our approach considerably faster.

We remark that Eigtool uses
the spectral norm $\|\cdot\|$ in \eqref{pseu}, while we apply the Frobenius norm for the 
matrices $E_i$. Since the $E_i$ are of rank one, they have the same spectral and Frobenius 
norms. 

\begin{figure}[ht]
\begin{center}
\includegraphics[width=6.5cm]{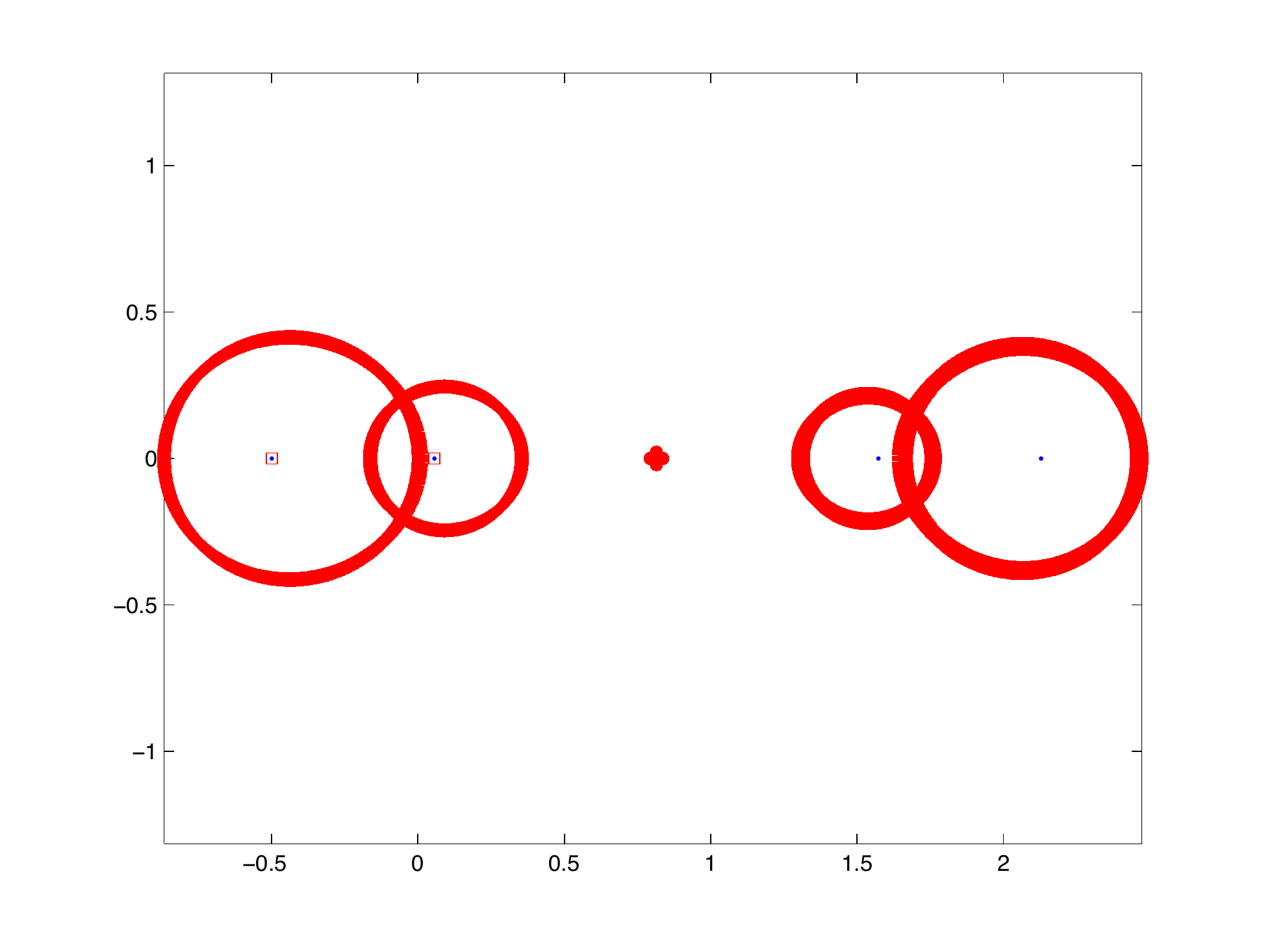}
\hskip 10mm
\includegraphics[width=6.5cm]{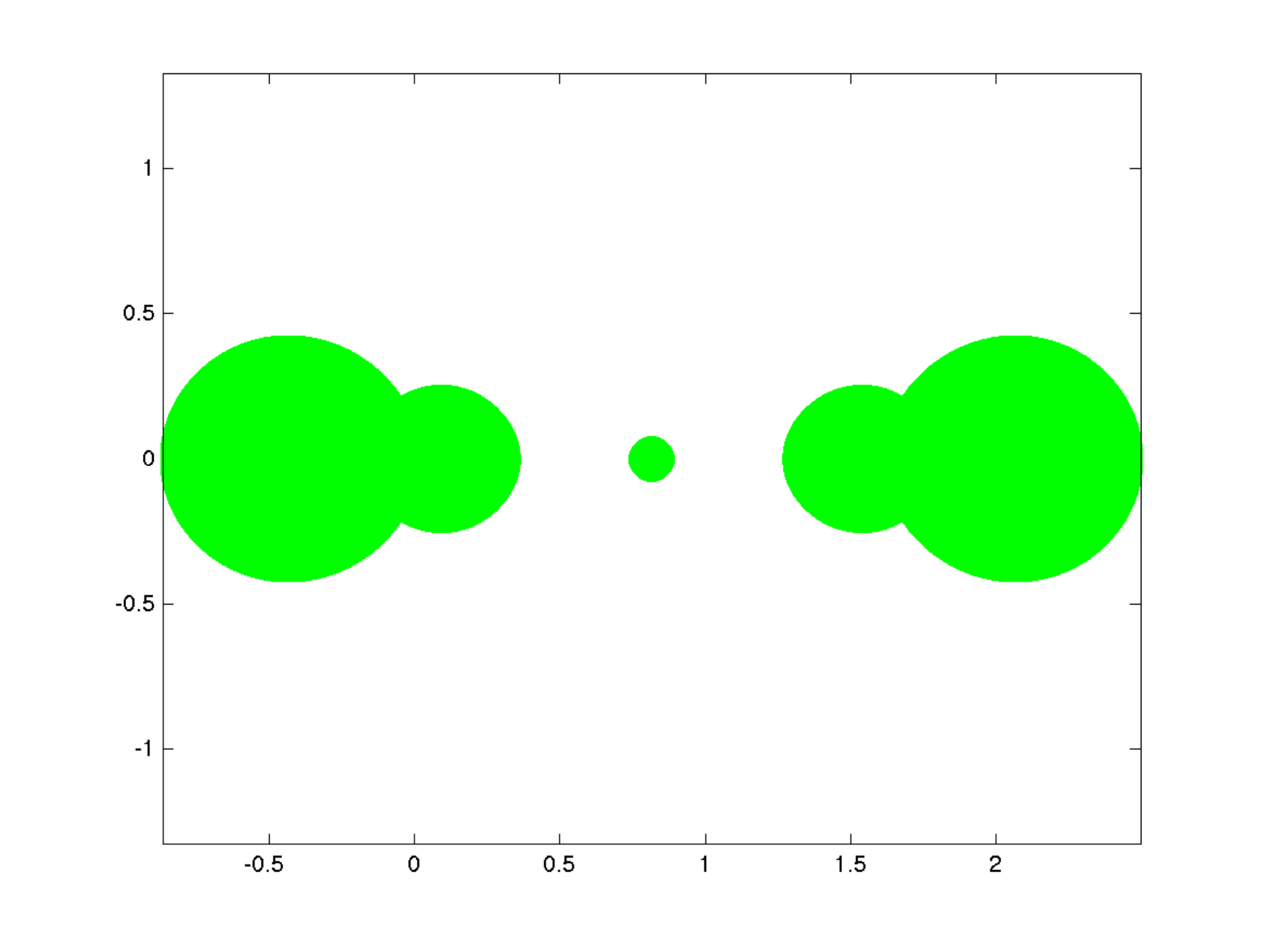}
\end{center}
\caption{Example 1. Left plot: $\Lambda^{\mc T}_{\varepsilon_2}(A)$ is approximated by the
eigenvalues of matrices of the form $A+\varepsilon_2 W^{\mc T}_1$ and 
$A+\varepsilon_2 W^{\mc T}_2$, where the $W^{\mc T}_j=W_j|_{\widehat{\mc T}}$ are 
normalized projected Wilkinson perturbations onto ${\mc T}$ associated with the 
eigenvalues $\lambda_j$, $j=1,2$ (marked by red squares), for $\eta:=\rme^{\rmi\theta_k}$ 
and $\theta_k:=2\pi(k-1)/10^3$, $k=1:10^3$, and $\varepsilon_2=10^{-0.8}$. Right plot:
$\Lambda^{\mc T}_{\varepsilon_2}(A)$  is approximated by the eigenvalues of matrices of 
the form $A+\varepsilon_2\rme^{\rmi\theta_k} E^{\mc T}_i$, $i,k=1:10^3$, where the 
$E^{\mc T}_i$ are unit-norm projected random perturbations in ${\mc T}$.\hfill\break}\label{fig4}
\end{figure} 

Next we turn to structured pseudospectra and perturbations. 
Let ${\mc T}$ be the space of tridiagonal Toeplitz matrices of order $5$.
We obtain from (\ref{rout_2}) the estimate 
$\varepsilon_2=10^{-0.8}$ of the structured distance from defectivity
$\varepsilon^{\mc T}_*$. It is achieved for the eigenvalues $\lambda_1$ and $\lambda_2$.
Figure \ref{fig4} (left) displays the spectra of matrices of the form
$A+\varepsilon_2 W^{\mc T}_1$ and $A+\varepsilon_2 W^{\mc T}_2$, where 
$W^{\mc T}_1=W_1|_{\widehat{\mc T}}$ and $W^{\mc T}_2=W_2|_{\widehat{\mc T}}$ are
normalized projected Wilkinson perturbations onto ${\mc T}$ associated with the 
eigenvalues $\lambda_1$ and $\lambda_2$, respectively, for $\eta:=\rme^{\rmi\theta_k}$ and
$\theta_k:=2\pi(k-1)/10^3$, $k=1:10^3$. The computations are described by Algorithm 2
with ${\cal S}={\cal T}$.

\begin{algorithm}
\DontPrintSemicolon
\KwData{matrix $A$, eigensystem $\{\lambda_i, x_i, y_i,\; \forall i=1:n\}$}
\KwResult{approximated $\Lambda_{\varepsilon^{\mc S}}^{\mc S}(A)$}

\nl compute $\varepsilon^{\mc S}$, $\{\hat\imath,\hat\jmath\}$   by \eqref{rout_2}\;
\nl compute  $W^{\mc S}_{\hat\imath}:=\varepsilon^{\mc S} \,y_{\hat\imath}x_{\hat\imath}^H|_{\widehat{\mc S}}$, $W^{\mc S}_{\hat\jmath}:=\varepsilon^{\mc S} \,y_{\hat\jmath}x_{\hat\jmath}^H|_{\widehat{\mc S}}$  \;
\nl display the spectrum of  $A+\eta W^{\mc S}_{\hat\imath}$ for $\eta:=\rme^{\rmi\theta_k}$, $\theta_k=2\pi(k-1)/10^3$, $k=1:10^3$ \;
\nl display the spectrum of  $A+\eta W^{\mc S}_{\hat\jmath}$ for $\eta:=\rme^{\rmi\theta_k}$, $\theta_k=2\pi(k-1)/10^3$, $k=1:10^3$ \;
\caption{Algorithm for computing an approximated structured pseudospectrum} \label{algo2}
\end{algorithm}

Figure \ref{fig4} (right) displays an
approximation of $\Lambda^{\mc T}_{\varepsilon_2}(A)$ given by the spectra of the matrices
$A+\varepsilon_2 \rme^{\rmi\theta_k} E^{\mc T}_i$, $i,k=1:10^3$, where the $E^{\mc T}_i$ 
are  random tridiagonal Toeplitz matrices scaled so that $\|E^{\mc T}_i\|_F=1$.
Eigtool \cite{Wr02} cannot be applied to determine structured pseudospectra. Notice that 
one component of 
the most $\Lambda_{\varepsilon_2}^{\mc T}$-sensitive pair of eigenvalues, $\lambda_2$, does 
not have one of the two largest structured condition numbers; see Table \ref{Tab1}. 
$\blacksquare$

\begin{table}[htb!]
\centering
\begin{tabular}{cccc}\hline
$i$ &$\lambda_i$ & $\kappa(\lambda_i)$ & $\kappa^{\mc T}(\lambda_i)$  \\ 
\hline 
${\phantom 0}1$ &${\phantom 0}5.4616 + 6.5356\,\rmi$ & $1.039\cdot 10^{1}$ & ${\phantom 0}1.169\cdot 10^{-1}$ \\
${\phantom 0}2$ &${\phantom 0}3.8552+ 5.1268\,\rmi$ & $2.999\cdot 10^{1}$ & ${\phantom 0}8.646\cdot 10^{-1}$ \\
${\phantom 0}3$ &${\phantom 0}1.7072 + 3.1264\,\rmi$ & $5.643\cdot 10^{1}$ & ${\phantom 0}5.665\cdot 10^{-1}$ \\
${\phantom 0}4$ &$-3.9451 - 0.1224\,\rmi$ & $1.534\cdot 10^{1}$ & $1.250\cdot 10^{0}$ \\
${\phantom 0}5$ &$-0.7339 - 3.2688\,\rmi$ & $2.528\cdot 10^{0}$ & ${\phantom 0}8.553\cdot 10^{-1}$ \\
${\phantom 0}6$ &${\phantom 0}0.3809 - 2.2234\,\rmi$ & $4.908\cdot 10^{0}$ & ${\phantom 0}6.596\cdot 10^{-1}$ \\
${\phantom 0}7$ &${\phantom 0}2.4409 - 0.7300\,\rmi$ & $2.373\cdot 10^{0}$ & ${\phantom 0}8.623\cdot 10^{-1}$ \\
${\phantom 0}8$ &${\phantom 0}1.5110- 1.0247\,\rmi$ & $8.071\cdot 10^{0}$ & ${\phantom 0}7.491\cdot 10^{-1}$ \\
${\phantom 0}9$ &$-2.2354 + 0.4417\,\rmi$ & $5.207\cdot 10^{1}$ & ${\phantom 0}9.748\cdot 10^{-1}$ \\
$10$ &$-0.2952 + 1.1966\,\rmi$ & $7.775\cdot 10^{1}$ & ${\phantom 0}3.750\cdot 10^{-1}$ \\
\hline
\end{tabular}
\caption{Example 2: Eigenvalue condition numbers.}
\label{Tab3}
\end{table}

\begin{figure}[ht]
\begin{center}
\includegraphics[width=6.5cm]{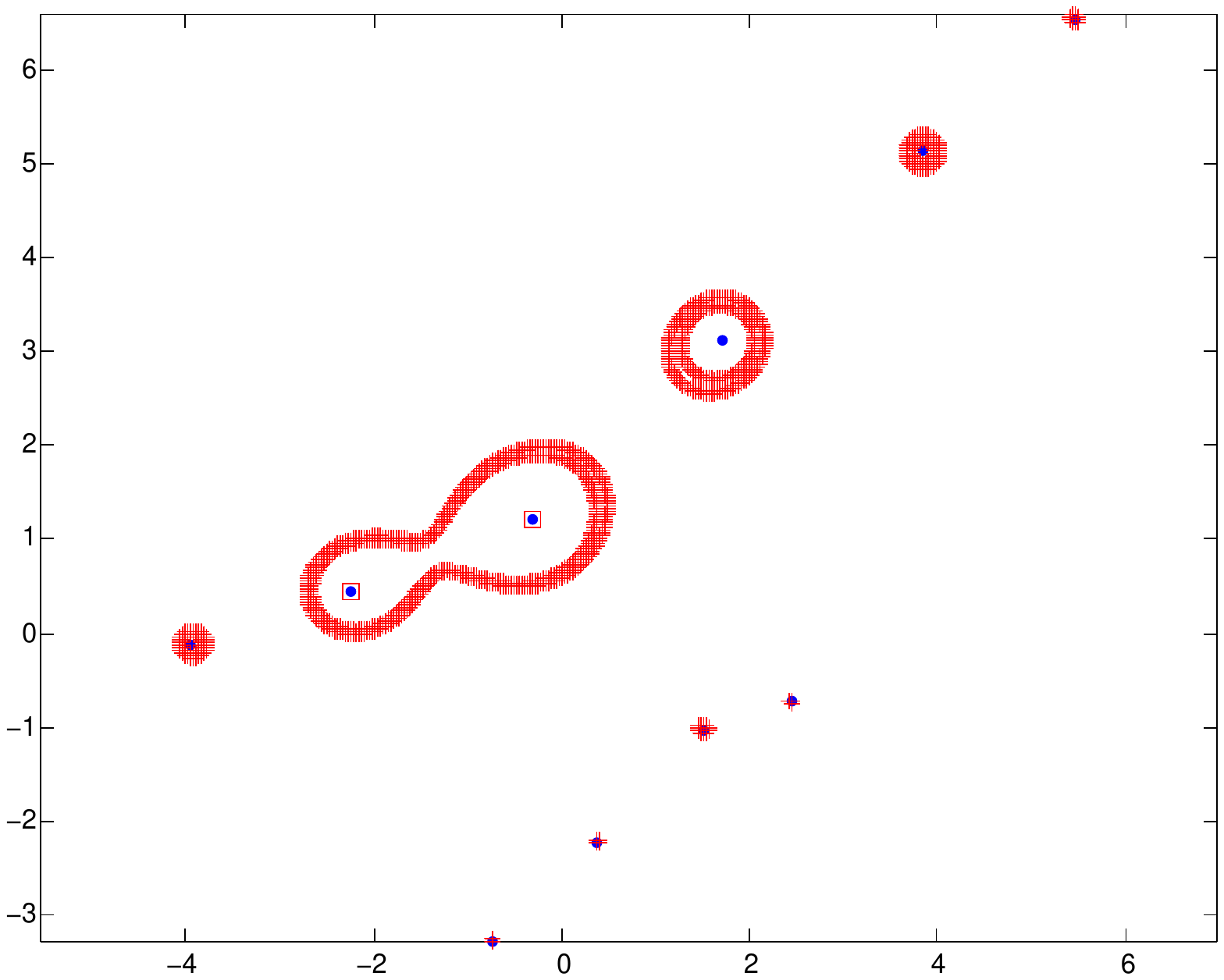}
\hskip 10mm
\includegraphics[width=6.5cm]{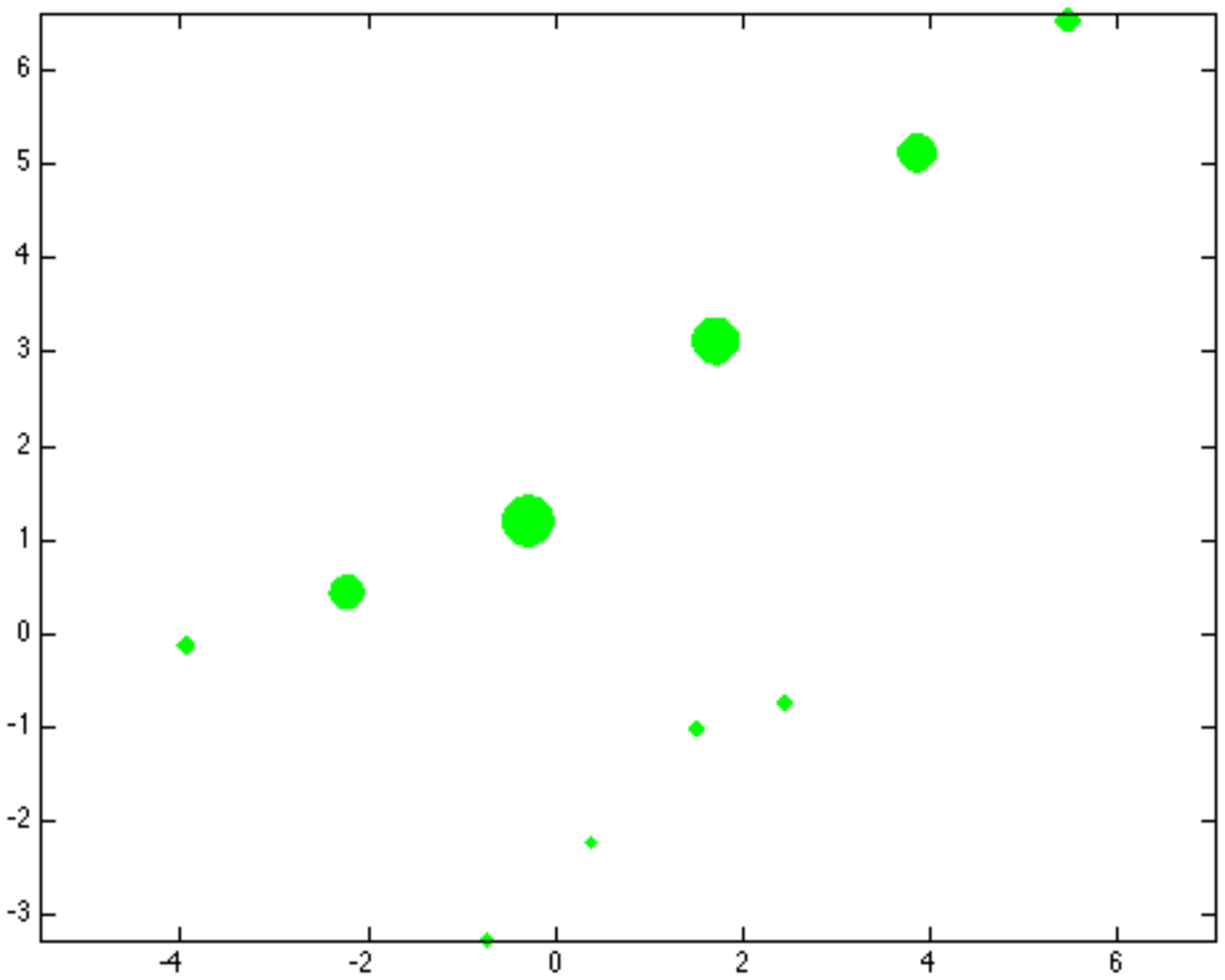}
\end{center}
\caption{Example 2. Left plot: $\Lambda_{\varepsilon_1}(A)$ is approximated by the 
eigenvalues of matrices of the form $A+\varepsilon_1 W_9$ and $A+\varepsilon_1 W_{10}$,
where the $W_j$ are Wilkinson perturbations associated with the eigenvalues $\lambda_j$,
$j=9,10$ (marked by red squares), for $\eta:=\rme^{\rmi\theta_k}$, 
$\theta_k:=2\pi(k-1)/10^3$, $k=1:10^3$, and $\varepsilon_1=10^{-2}$. Right plot: 
$\Lambda_{\varepsilon_1}(A)$  is approximated by the eigenvalues of matrices of the form 
$A+\varepsilon_1\rme^{\rmi\theta_k} E_i$, $i,k=1:10^3$, where the $E_i$ are unit-norm 
rank-one random perturbations.\hfill\break}\label{fig1_p}
\end{figure}

\begin{figure}[ht]
\begin{center}
\includegraphics[width=13cm]{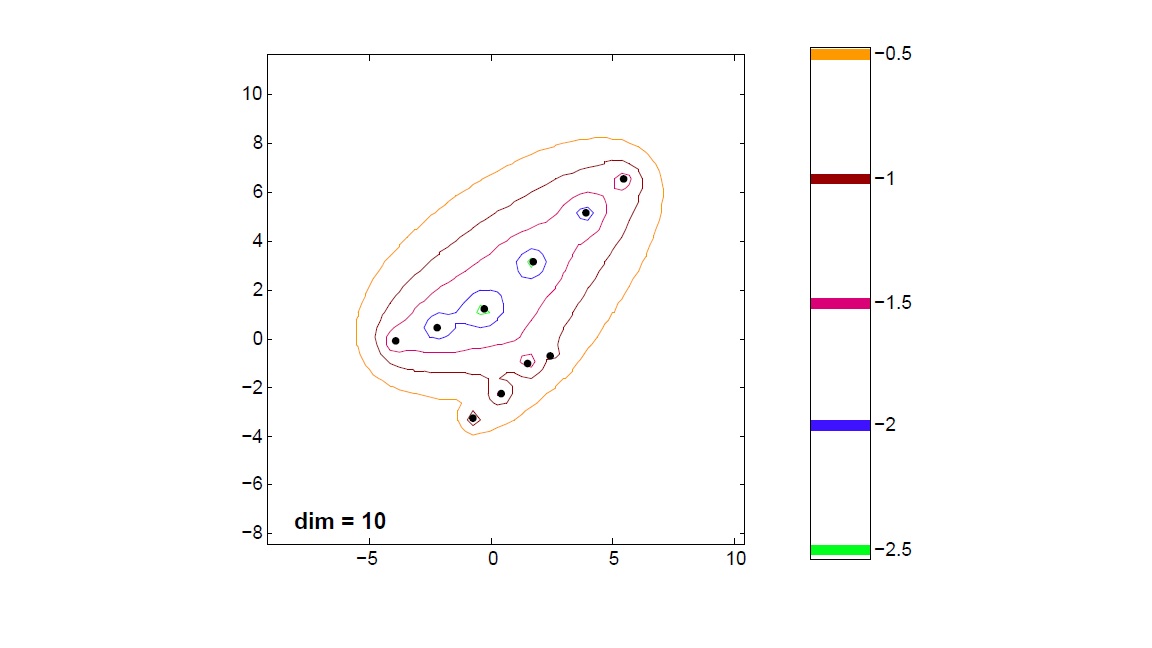}
\caption{Example 2: $\varepsilon$-pseudospectra by Eigtool, where $\varepsilon=10^k$, 
$k=-2.5:0.5:-0.5$. dim shows the order of the matrix.\hfill\break}\label{fig3_p}
\end{center}
\end{figure}

{\bf Example 2}. We consider a complex pentadiagonal Toeplitz matrix of order $n=10$ 
constructed analogously as the matrix of Example 1. Traditional and structured condition
numbers for the eigenvalues are shown in Table \ref{Tab3}. The estimate (\ref{rout_1}) of
the distance to defectivity is $\varepsilon_1=10^{-2}$; it is achieved for the indices $9$ 
and $10$. Figure \ref{fig1_p} (left) displays the spectra of matrices of the form
$A+\varepsilon_1 W_9$ and $A+\varepsilon_1 W_{10}$, where the $W_j$ are Wilkinson 
perturbations \eqref{rank1} associated with the eigenvalues $\lambda_j$, $j=9,10$, for
$\eta:=\rme^{\rmi\theta_k}$ and $\theta_k:=2\pi(k-1)/10^3$, $k=1:10^3$. Figure 
\ref{fig1_p} (right) displays the approximated $\varepsilon_1$-pseudospectrum given by the 
spectra of the matrices $A+\varepsilon_1\rme^{\rmi\theta_k}E_i$, $i,k=1:10^3$, 
where the $E_i$ are random  rank-one perturbations with $\|E_i\|_F=1$. 
Figure \ref{fig3_p} 
depicts pseudospectra determined by Eigtool \cite{Wr02}. Comparing the
$\varepsilon_1$-pseudospectrum of Figure \ref{fig3_p} with Figure \ref{fig1_p} shows that 
the simple computations of this paper, based on Wilkinson perturbations \eqref{rank1} and
illustrated by Figure \ref{fig1_p} (left), can give more accurate approximations of 
pseudospectra and require less computational effort than the approach used for Figure
\ref{fig1_p} (right). Notice that the most $\Lambda_{\varepsilon_1}$-sensitive pair of 
eigenvalues do not have the largest (unstructured) condition numbers; see Table \ref{Tab3}.

\begin{figure}[ht]
\begin{center}
\includegraphics[width=6.5cm]{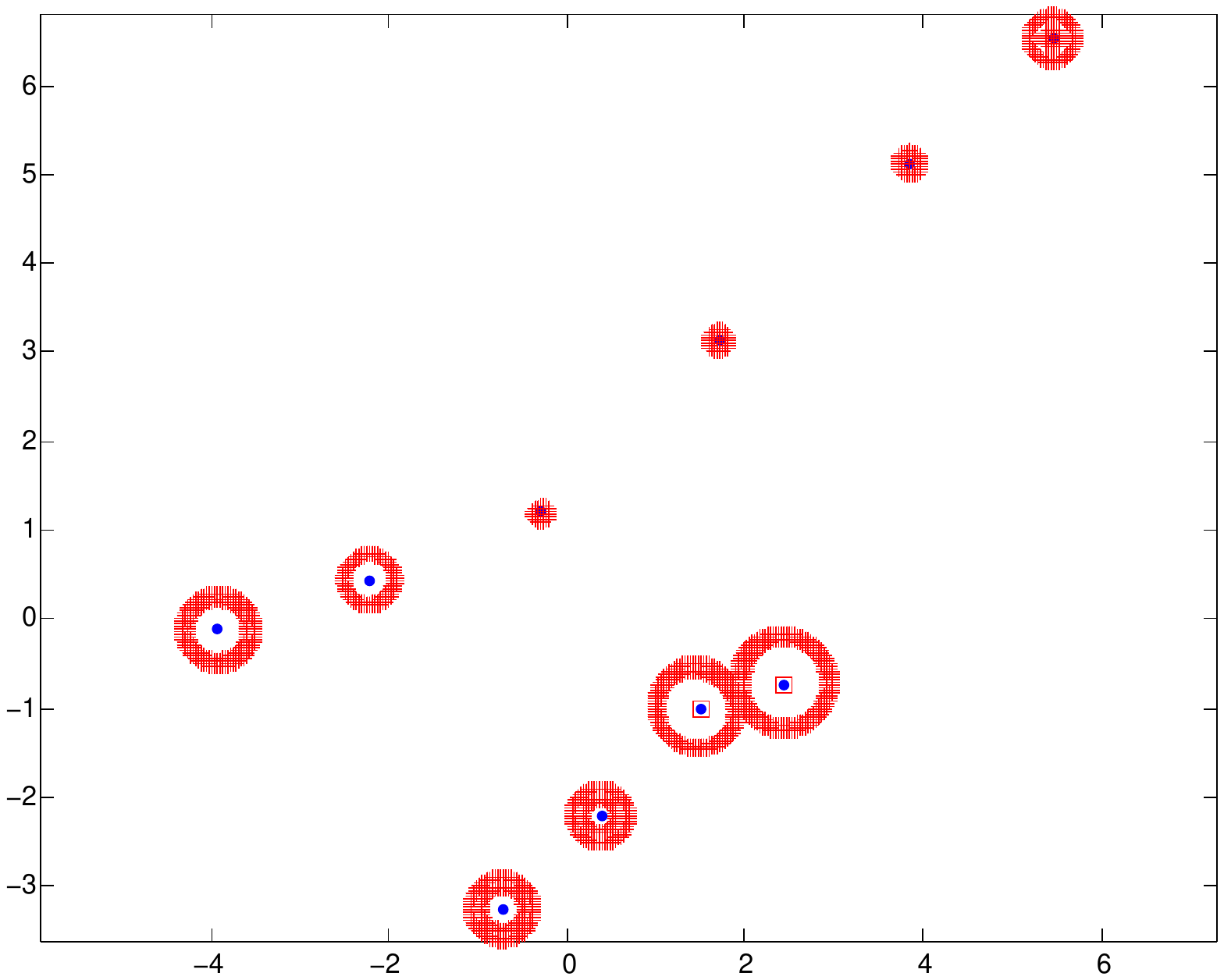}
\hskip 10mm
\includegraphics[width=6.5cm]{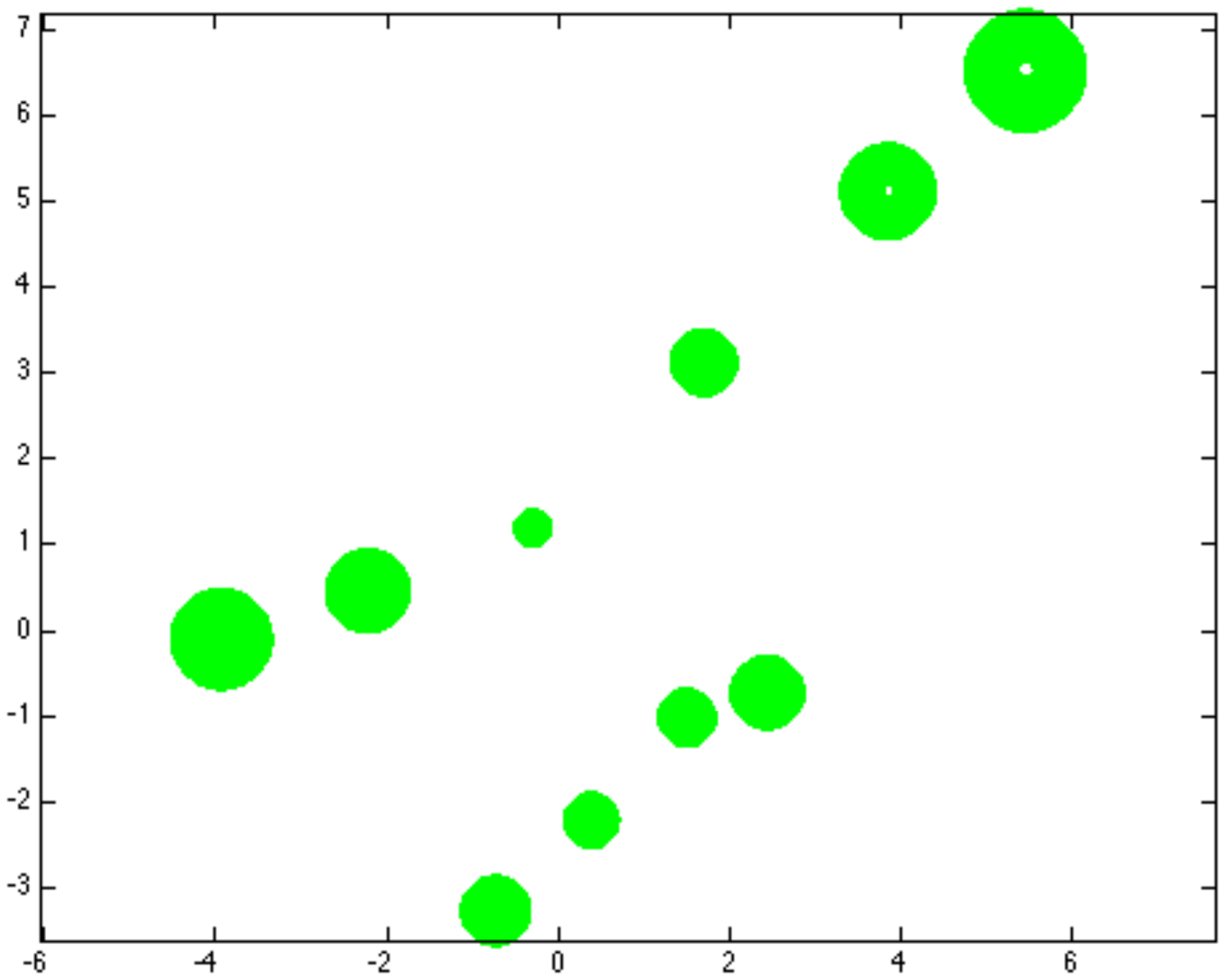}
\end{center}
\caption{Example 2. Left plot: $\Lambda^{\mc T}_{\varepsilon_2}(A)$ is approximated by the
eigenvalues of matrices of the form $A+\varepsilon_2 W^{\mc T}_7$ and 
$A+\varepsilon_2 W^{\mc T}_8$, where the $W^{\mc T}_j=W_j|_{\widehat{\mc T}}$ are 
normalized projected Wilkinson perturbations onto ${\mc T}$ associated with the 
eigenvalues $\lambda_j$, $j=7,8$ (marked by red squares), for $\eta:=\rme^{\rmi\theta_k}$
and $\theta_k:=2\pi(k-1)/10^3$, $k=1:10^3$, and $\varepsilon_2=10^{-0.2}$. Right plot: 
$\Lambda^{\mc T}_{\varepsilon_2}(A)$  is approximated by the eigenvalues of matrices of the
form $A+\varepsilon_2 \rme^{\rmi\theta_k} E^{\mc T}_i$, $i,k=1:10^3$, where the 
$E^{\mc T}_i$ are unit-norm pentadiagonal Toeplitz random perturbations.\hfill\break}
\label{fig2_p}
\end{figure}

We now consider structured perturbations and pseudospectra. Let ${\mc T}$ be the space 
of pentadiagonal Toeplitz matrices of order $10$. We obtain from (\ref{rout_2})
the estimate $\varepsilon_2=10^{-0.2}$ of the structured distance from defectivity 
$\varepsilon^{\mc T}_*$. It is achieved for the eigenvalues $\lambda_7$ and $\lambda_8$. 
Figure \ref{fig2_p} (left) displays the spectra of matrices of the form
$A+\varepsilon_2 W^{\mc T}_7$ and $A+\varepsilon_2 W^{\mc T}_8$, where the 
$W^{\mc T}_j=W_j|_{\widehat{\mc T}}$ are normalized projected Wilkinson perturbations onto
${\mc T}$ associated with the eigenvalues $\lambda_j$ for $j=7,8$ with 
$\eta:=\rme^{\rmi\theta_k}$ for $\theta_k=2\pi(k-1)/10^3$, $k=1:10^3$.
Figure \ref{fig2_p}
(right) displays the approximation of $\Lambda^{\mc T}_{\varepsilon_2}(A)$ given by the 
spectra of the matrices $A+\varepsilon_2 \rme^{\rmi\theta_k} E^{\mc T}_i$, $i,k=1:10^3$, 
where the $E^{\mc T}_i$ are random  pentadiagonal Toeplitz matrices scaled so that 
$\|E^{\mc T}_i\|_F=1$. According to Table \ref{Tab3}, the most 
$\Lambda_{\varepsilon_2}^{\mc T}$-sensitive pair of eigenvalues do not have the largest 
structured condition numbers.  $\blacksquare$

\begin{table}[htb!]
\centering
\begin{tabular}{cccc}\hline
$i$ &$\lambda_i$ & $\kappa(\lambda_i)$ & $\kappa^{\mc H}(\lambda_i)$  \\ 
\hline 
$1$&$-2.0595$ & $1.092\cdot 10^{0}$ & ${\phantom 0}7.725\cdot 10^{-1}$\\
$2$&${\phantom 0}2.0595$ & $1.092\cdot 10^{0}$ & ${\phantom 0}7.725\cdot 10^{-1}$ \\
$3$&$-0.6686$ & $1.758\cdot 10^{0}$ & $1.252\cdot 10^{0}$ \\
$4$&${\phantom 0}0.6686$ & $1.758\cdot 10^{0}$ & $1.252\cdot 10^{0}$ \\
$5$&${\phantom 0}0.3677$ & $4.097\cdot 10^{0}$ & $2.926\cdot 10^{0}$ \\
$6$&${\phantom 0}0.2151$ & $3.958\cdot 10^{0}$ & $3.009\cdot 10^{0}$ \\
$7$&$-0.2151$ & $3.958\cdot 10^{0}$ & $3.009\cdot 10^{0}$ \\
$8$&$-0.3677$ & $4.097\cdot 10^{0}$ & $2.926\cdot 10^{0}$ \\
\hline
\end{tabular}
\caption{Example 3: Eigenvalue condition numbers.}
\label{Tab4}
\end{table}

\begin{figure}[ht]
\begin{center}
\includegraphics[width=6.5cm]{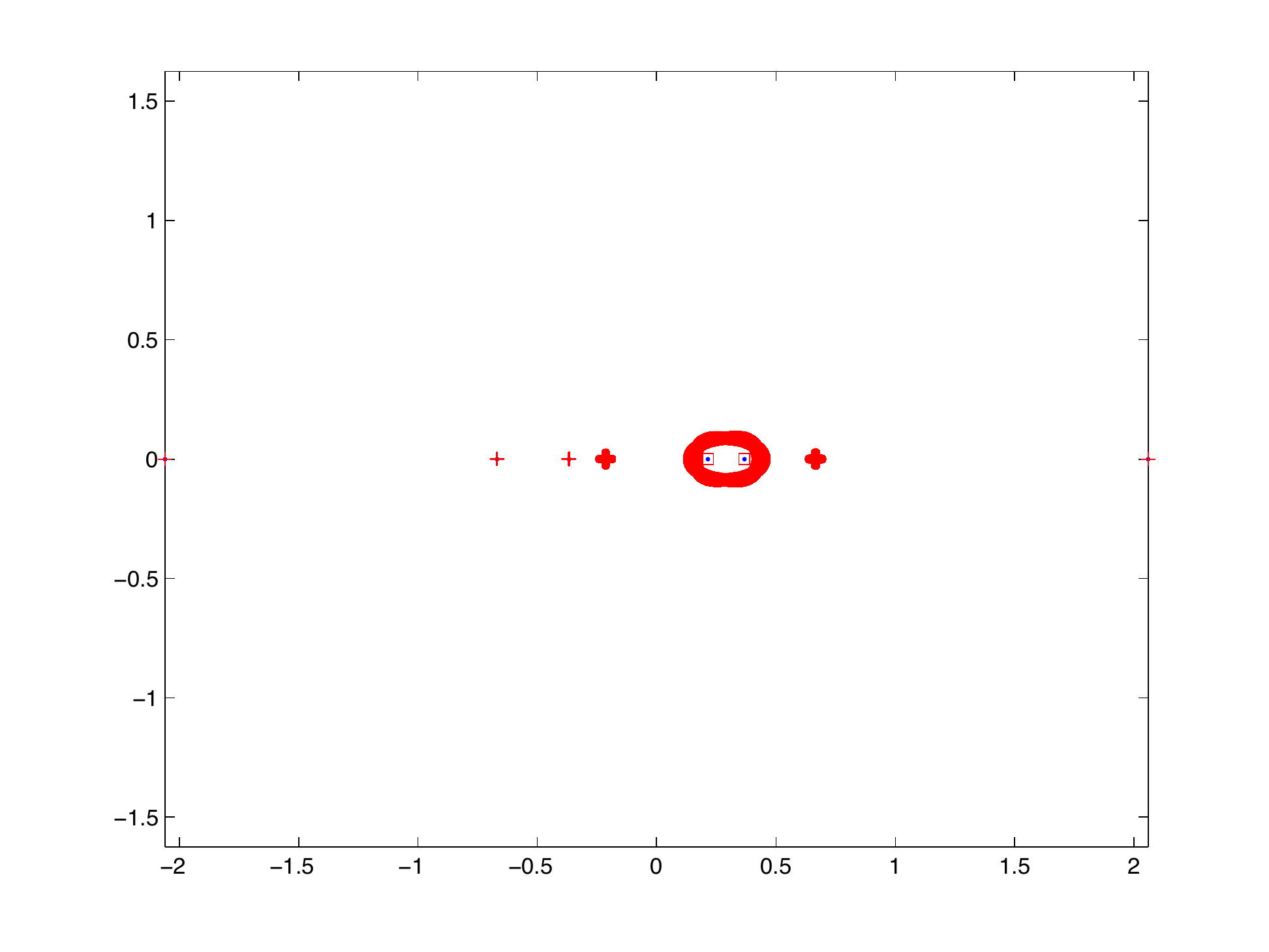}
\hskip 10mm
\includegraphics[width=6.5cm]{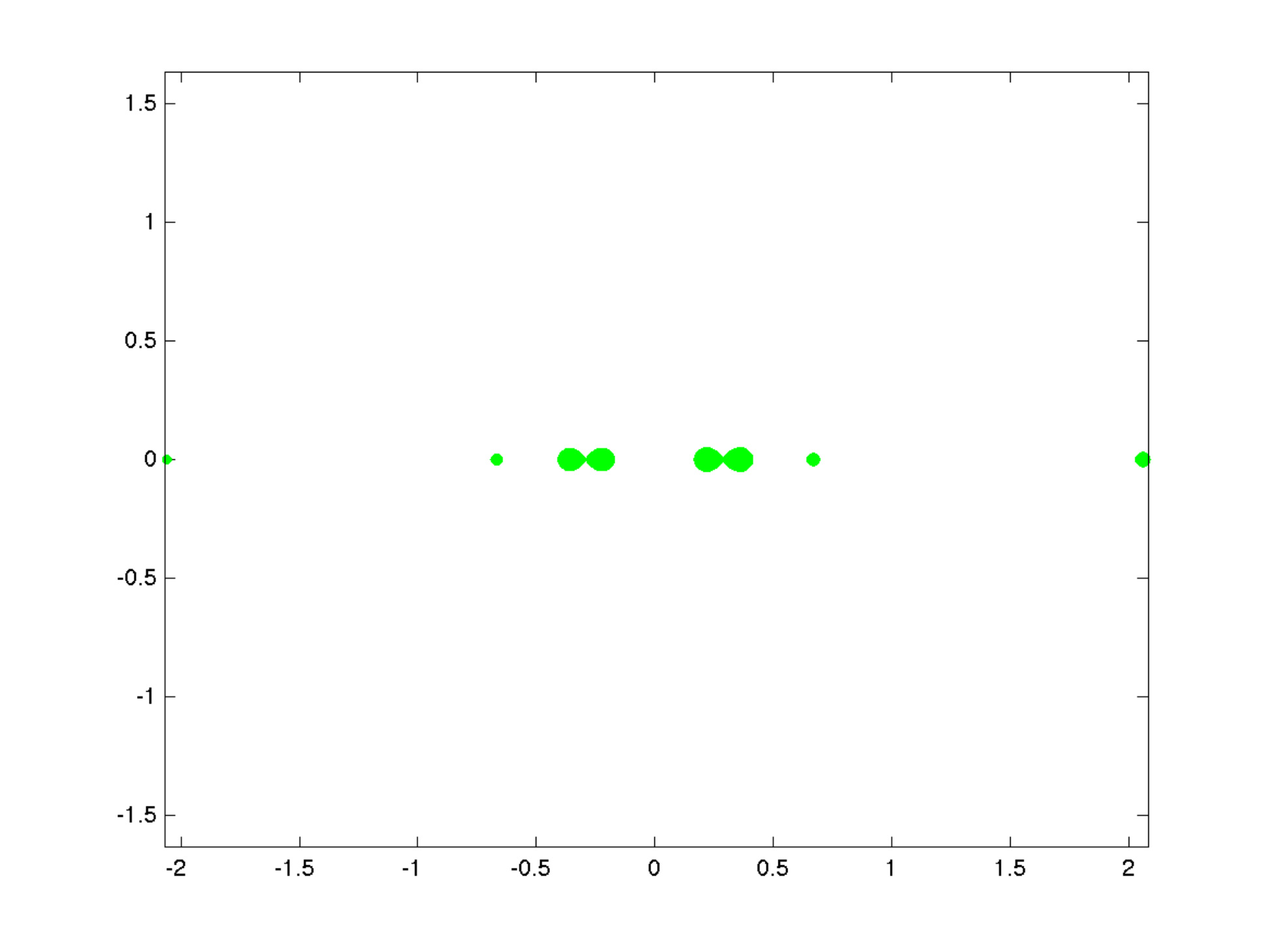}
\end{center}
\caption{Example 3. Left plot: $\Lambda_{\varepsilon_1}(A)$ is approximated by the 
eigenvalues of matrices of the form $A+\varepsilon_1 W_5$ and $A+\varepsilon_1 W_6$, where
the $W_j$ are Wilkinson perturbations associated with the eigenvalues $\lambda_j$, $j=5,6$
(marked by red squares), for $\eta:=\rme^{\rmi\theta_k}$, 
$\theta_k:=2\pi(k-1)/10^3$, $k=1:10^3$, and $\varepsilon_1=10^{-1.6}$. Right plot: 
$\Lambda_{\varepsilon_1}(A)$  is approximated by the eigenvalues of 
$A+\varepsilon_1\rme^{\rmi\theta_k}E_i$, $i,k=1:10^3$, where the $E_i$ are unit-norm
rank-one random perturbations.\hfill\break}\label{fig1_3}
\end{figure}

\begin{figure}[ht]
\begin{center}
\includegraphics[width=13cm]{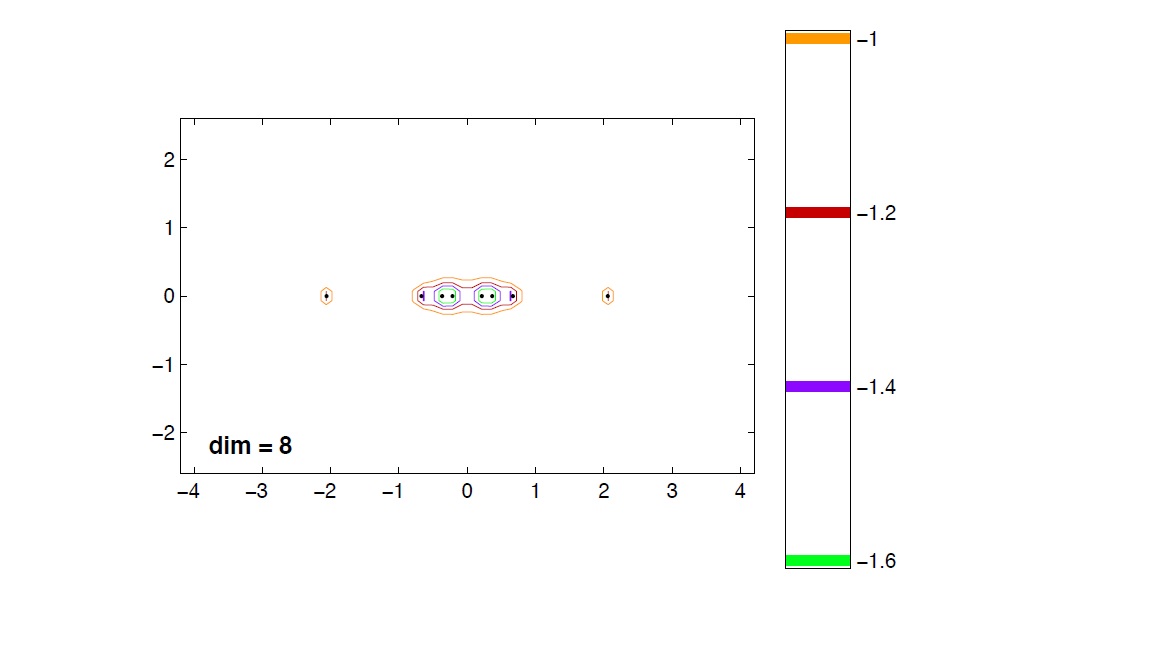}
\caption{Example 3: $\varepsilon$-pseudospectra by Eigtool, where $\varepsilon=10^k$, 
$k=-1.6:0.2:-1$. dim shows the order of the matrix.\hfill\break}\label{fig2_3}
\end{center}
\end{figure}

\subsection{Hamiltonian structure} 
$\\$ {\bf Example 3}. Let $A=M|_{{\mc H}}$, where $M\in\mathbb{R}^{8\times 8}$ has random 
entries, i.e., $A$ is the closest Hamiltonian matrix to $M$; cf. Proposition \ref{prop:1}.
The traditional and structured condition numbers for the eigenvalues of $A$ are shown in
Table \ref{Tab4}. We obtain from (\ref{rout_1}) the upper bound $\varepsilon_1=10^{-1.6}$
for the distance to defectivity $\varepsilon_*$. It is achieved for the indices $5$ and 
$6$. Figure \ref{fig1_3} (left) is analogous to Figure \ref{fig1_p}; it displays the 
spectra of matrices of the form $A+\varepsilon_1 W_5$ and $A+\varepsilon_1 W_6$, where the
$W_j$ are Wilkinson perturbations \eqref{rank1} associated with the eigenvalues 
$\lambda_j$, $j=5,6$, for $\eta:=\rme^{\rmi\theta_k}$, $\theta_k:=2\pi(k-1)/10^3$, 
$k=1:10^3$. Figure \ref{fig1_3} (right) shows the approximated 
$\varepsilon_1$-pseudospectrum given by the spectra of the matrices
$A+\varepsilon_1 \rme^{\rmi\theta_k}E_i$, $i,k=1:10^3$, where the $E_i$ are random
rank-one perturbations with $\|E_i\|_F=1$. Figure \ref{fig2_3} displays graphs produced by
Eigtool \cite{Wr02}. A comparison of Figures \ref{fig1_3} and \ref{fig2_3} shows the 
effectiveness of using (\ref{rout_1}) to identify pertinent eigenvalue pairs. Similarly as 
in the preceding examples, $2\cdot 10^3$ Wilkinson perturbations yield a much better 
approximation of the $\varepsilon_1$-pseudospectrum than a simulation with $1\cdot 10^6$ 
Hamiltonian random perturbations. The latter simulation does not show coalescence of 
components of the $\varepsilon_1$-pseudospectrum.

\begin{figure}[ht]
\begin{center}
\includegraphics[width=6.5cm]{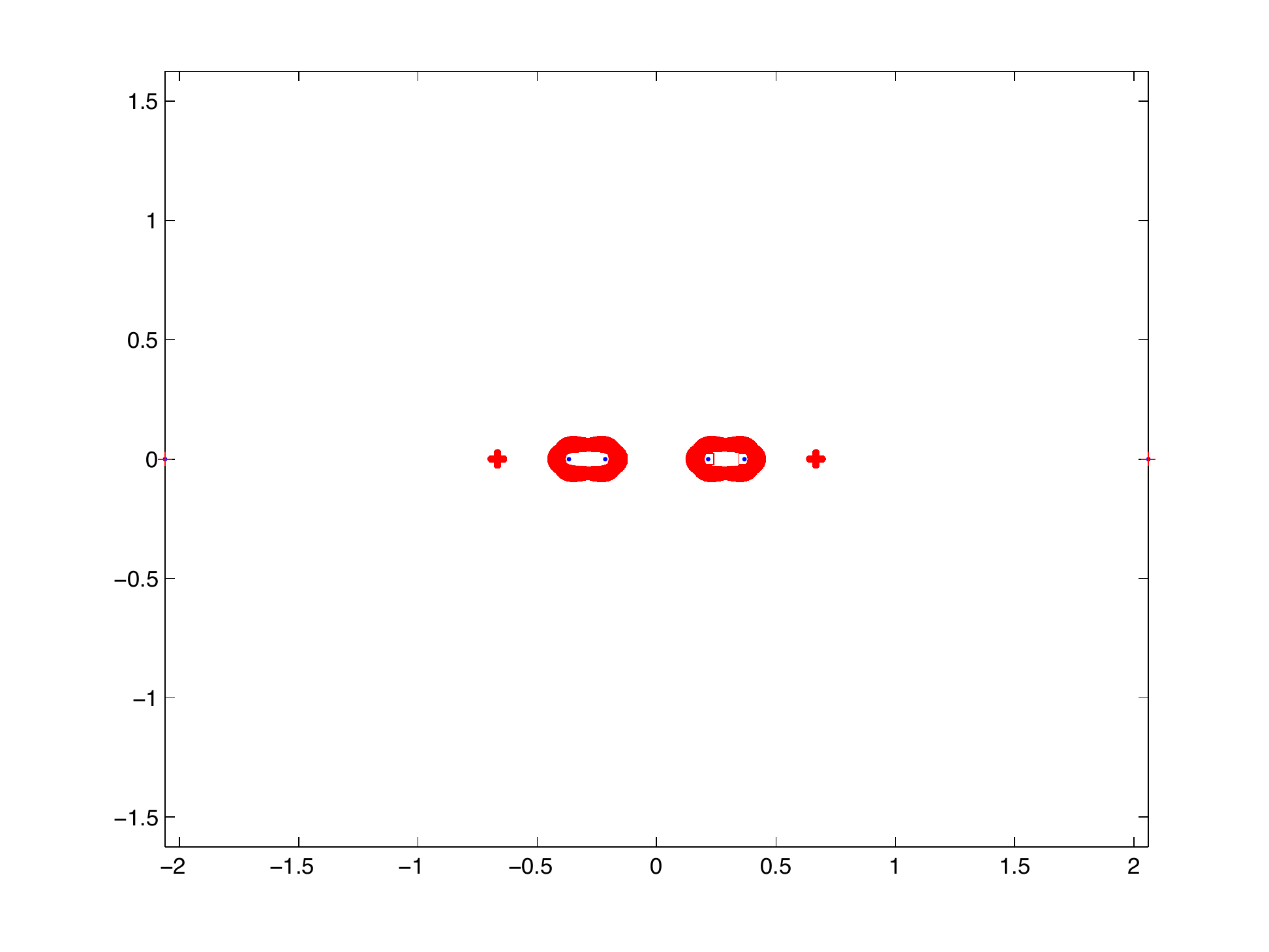}
\hskip 10mm
\includegraphics[width=6.5cm]{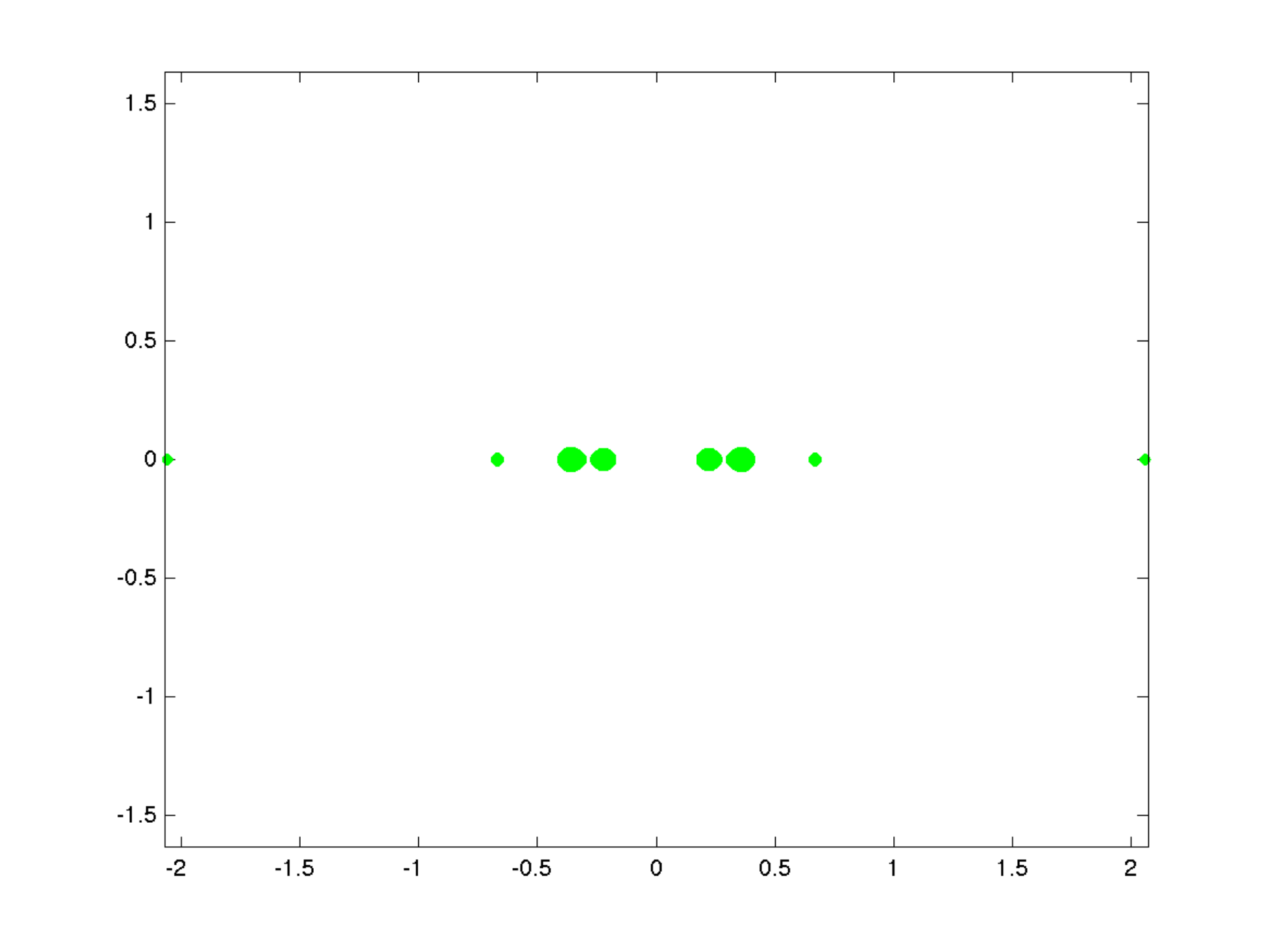}
\end{center}
\caption{Example 3. Left plot: $\Lambda^{\mc H}_{\varepsilon_2}(A)$ is approximated by the
eigenvalues of matrices of the form  $A+\varepsilon_2 W^{\mc H}_5$ and 
$A+\varepsilon_2 W^{\mc H}_6$, where the $W^{\mc T}_j=W_j|_{\widehat{\mc T}}$ are 
normalized projected Wilkinson perturbations onto ${\mc H}$ associated with the 
eigenvalues $\lambda_j$, $j=5,6$ (marked by red squares), for $\eta:=\rme^{\rmi\theta_k}$,
$\theta_k:=2\pi(k-1)/10^3$, $k=1:10^3$, and $\varepsilon_2=10^{-1.6}$. Right plot: 
$\Lambda^{\mc H}_{\varepsilon_2}(A)$ is approximated by the eigenvalues of 
$A+\varepsilon_2\rme^{\rmi\theta_k} E^{\mc H}_i$, $i,k=1:10^3$, where the $E^{\mc H}_i$ 
are unit-norm Hamiltonian random perturbations.\hfill\break}\label{fig3_3}
\end{figure}

We now consider structured perturbations and evaluate (\ref{rout_2}). This gives 
$\varepsilon_2=10^{-1.6}$, which is the same as $\varepsilon_1$ above. We obtain 
the same indices, $5$ and $6$, as for the unstructured situation. Figure \ref{fig3_3} is
analogous to Figure \ref{fig1_p}. Thus, Figure \ref{fig3_3} (left) displays spectra of 
matrices of the form $A+\varepsilon_2 W^{\mc H}_5$ and $A+\varepsilon_2 W^{\mc H}_6$, 
where the $W^{\mc H}_j=W_j|_{\widehat{\mc H}}$ are normalized projected Wilkinson 
perturbations onto the space ${\mc H}$ of the real Hamiltonian matrices associated with
the eigenvalues $\lambda_j$, $j=5,6$, with $\eta:=\rme^{\rmi\theta_k}$, 
$\theta_k:=2\pi(k-1)/10^3$, $k=1:10^3$. 
It is shown in \cite{BN} that in this case, the real worst-case perturbations are
rank-2 matrices.
In fact, one has that just the perturbation $\varepsilon_2 W^{\mc H}_5$ with $\eta=-1$,
or even just the perturbation $\varepsilon_2 W^{\mc H}_6$ with $\eta=1$,  suffices to cause 
coalescence of the pairs $(\lambda_7,\lambda_8)$ and $(\lambda_5,\lambda_6)$; see (\ref{p:2}).
Figure \ref{fig3_3} (right) shows an approximation of the spectrum
$\Lambda_{\varepsilon_2}^{\mc H}(A)$ determined by the spectra of the matrices 
$A+\varepsilon_2 \rme^{\rmi\theta_k} E^{\mc H}_i$, $i,k=1:10^3$, where the
$E^{\mc H}_i$ are Hamiltonian random perturbations scaled so that $\|E^{\mc H}_i\|_F=1$.

We observe that for both the structured and unstructured cases, the use of (\ref{rout_1})
together with Wilkinson perturbations, or of (\ref{rout_2}) with projected Wilkinson
perturbations, give more accurate approximations of the 
$\varepsilon_1$-pseudospectrum and the structured $\varepsilon_2$-pseudospectrum, 
respectively, than using many more unstructured or structured random perturbations. Moreover,
the random perturbations do not provide information about coalescence of components of the
unstructured or structured pseudospectra.

Note that, although the most $\Lambda_{\varepsilon_1}$-sensitive pair of eigenvalues coincide with
the most $\Lambda_{\varepsilon_2}^{\mc H}$-sensitive pair of eigenvalues, see Table \ref{Tab4}, 
the structured approach has the advantage of preserving eigenvalue symmetries in finite 
precision arithmetic,
as it illustrated by left-hand side plots in Figures \ref{fig1_3} and \ref{fig3_3}.

\begin{figure}[ht]
\begin{center}
\includegraphics[width=6.5cm]{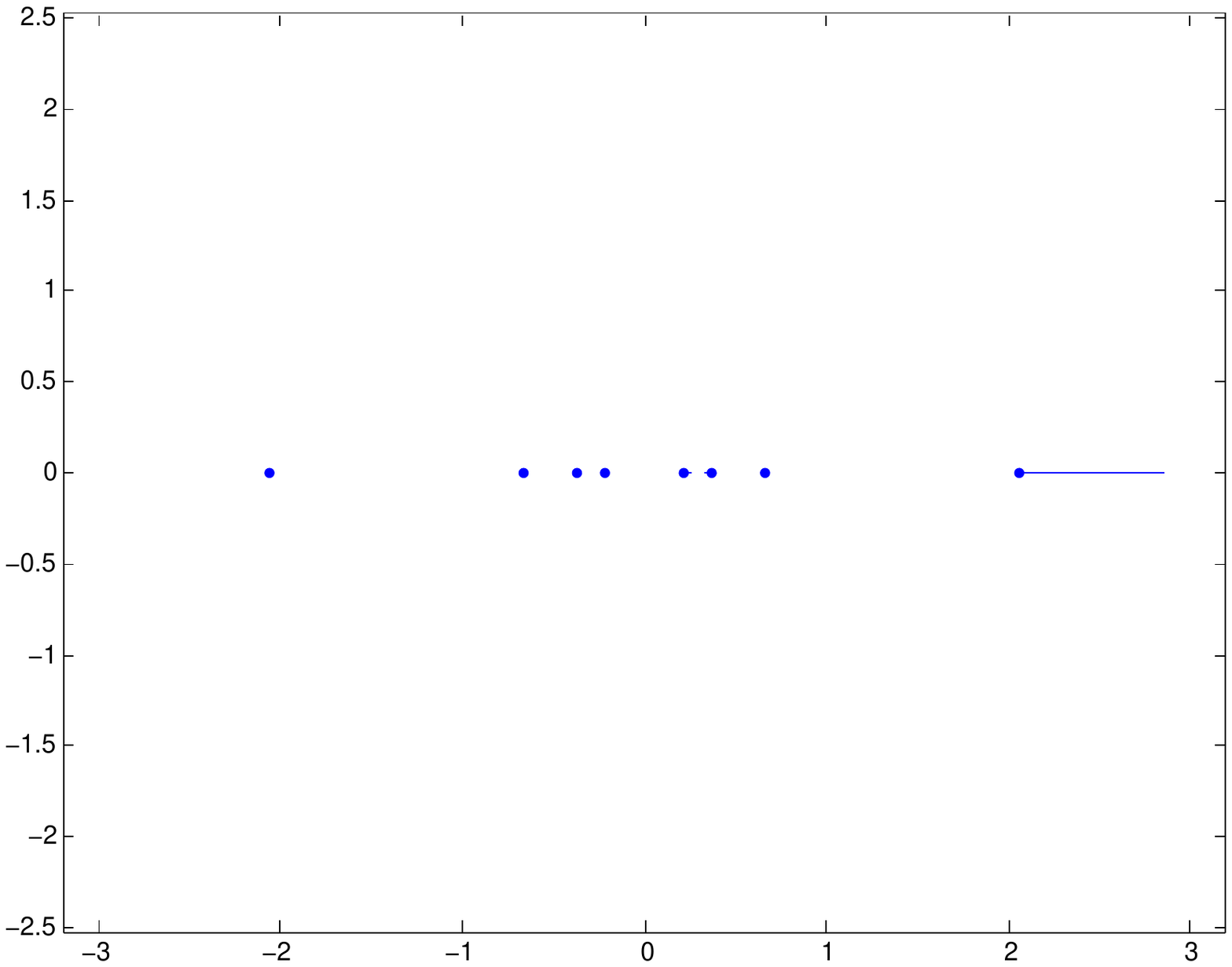}
\hskip 10mm
\includegraphics[width=6.5cm]{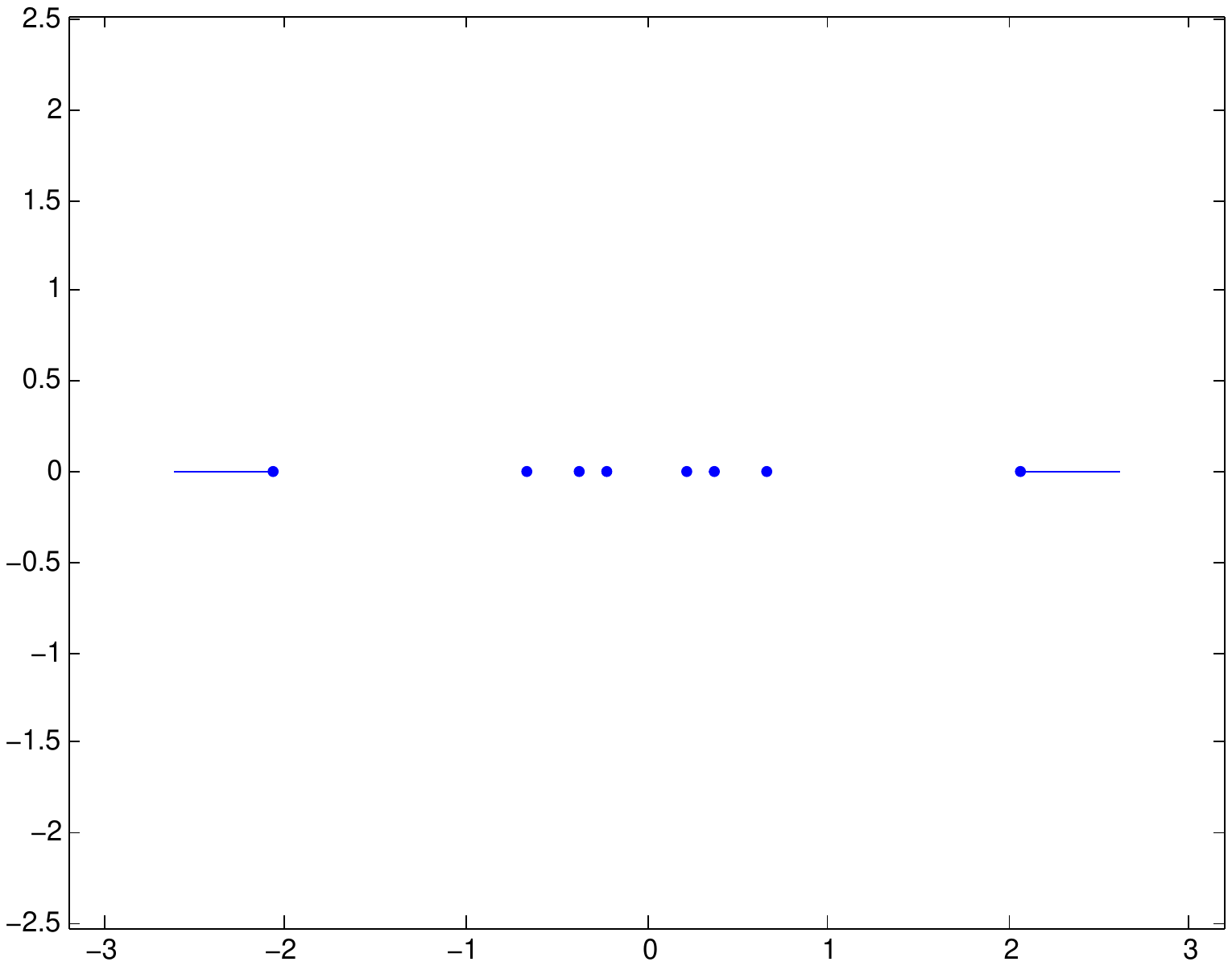}
\end{center}
\caption{Example 3: Approximated pseudo-eigenvalues and structured pseudo-eigenvalues. 
Left plot: $E$ has entries $1/n$ and $\varepsilon$ increases from $0$ to 
$10^{-0.1}$. Right plot: perturbation $E|_{\widehat{\mc H}}$ and the same $\varepsilon$.
\hfill\break}\label{fig13}
\end{figure}

We conclude this example with another illustration of the sensitivity of the eigenvalues 
of $A$ to perturbations. Let $E$ be a rank-one matrix of norm one. Figure \ref{fig13} 
shows for all eigenvalues $\lambda_i=\lambda_i(0)$, $1\leq i\leq 8$, of $A$, the behavior 
of the pseudo-eigenvalues $\lambda_i(\varepsilon)$ of $A+\varepsilon E$ and of the structured
pseudo-eigenvalues $\lambda_i^{\mc H}(\varepsilon)$ of $A+\varepsilon E^{\mc H}$, where
\[
\lambda_i(\varepsilon)\approx\lambda_i(0)+\varepsilon\frac{y_i^HEx_i}{y_i^Hx_i},\qquad
\lambda_i^{\mc H}(\varepsilon)\approx\lambda_i(0)+\varepsilon
\frac{y_i^HE^{\mc H}x_i}{y_i^Hx_i},
\]
and $y_i^Hx_i$ is real and positive. Here $x_i$ and $y_i$ are right and left unit 
eigenvectors associated with the eigenvalue $\lambda_i(0)$ and $\varepsilon$ increases 
from $0$ to $10^{-0.1}$. In Figure \ref{fig13} (left) the perturbation matrix $E$ is the 
rank-one perturbation with all elements equal to $1/8$, and in Figure 
\ref{fig13} (right) the structured perturbation matrix is 
$E^{\mc H}=E|_{\widehat{\mc H}}$. Real Hamiltonian matrices have eigenvalues
in $\pm$ pairs. This can be seen in Figure \ref{fig13} (right). Also, the rightmost
eigenvalue is perturbed less under the structured perturbation. 
Rough lower bounds for the $\varepsilon$-pseudospectral abscissa and for its structured version can be easily 
deduced.
$\blacksquare$

\section{Conclusions and remarks}\label{sec:5}
The computed examples illustrate that standard (unstructured) pseudospectra can be 
well approximated by using 
suitable 
``worst case'' rank-one perturbations of the given matrix, 
i.e., by using Wilkinson perturbations associated with the two eigenvalues whose 
pseudospectral components are likely to first coalesce, as determined by (\ref{rout_1}).
For structured matrices, such as banded non-Hermitian Toeplitz matrices or Hamiltonian 
matrices, the structured pseudospectra can be well approximated by using normalized 
projections of Wilkinson perturbations associated with two eigenvalues whose components
in the structured pseudospectra are likely to first coalesce, as determined by 
(\ref{rout_2}). For the Hamiltonian-structured case, this strategy gives rise to rank-two
approximated structured pseudospectra, since the Hamiltonian projection of a Wilkinson 
perturbation is of rank two. Finally, a simple strategy for approximating the structured
$\varepsilon$-pseudospectral abscissa [or radius]  (with respect to the Frobenius norm) consists of 
perturbing the original matrix $A$ by the $\varepsilon$-normalized projected Wilkinson 
perturbation associated 
with 
the rightmost [or largest] eigenvalue. We noticed in computations
that one often gets an extremely cheap, though quite rough, lower bound for the 
structured $\varepsilon$-pseudospectral abscissa by computing 
the real part of
the rightmost eigenvalue of $A$ perturbed by 
the $\varepsilon$-normalized projected all-ones matrix.

\section*{Acknowledgment}
The authors would like to thank the referees for suggestions that improved the presentation.

\end{document}